\documentclass[reqno, 11pt]{amsart}
\usepackage{amsmath}
\usepackage{amscd}
\usepackage{graphics}
\usepackage{latexsym}
\usepackage{hyperref}

\theoremstyle{plain}
\newtheorem{theorem}{Theorem}[section]
\newtheorem{thm}[theorem]{Theorem}

\newtheorem{lem}[theorem]{Lemma}
\newtheorem{prop}[theorem]{Proposition}

\theoremstyle{definition}
\newtheorem{defn}[theorem]{Definition}
\newtheorem{notat}[theorem]{Notation}

\newtheorem{hyp}[theorem]{Hypothesis}

\newtheorem{prin}[theorem]{Principle}

\newtheorem{rmk}[theorem]{Remark}

\theoremstyle{remark}

\newcommand{\ZZ}{\mathbb{Z}}
\newcommand{\CC}{\mathbb{C}}

\newcommand{\PP}{\mathbb{P}}
\newcommand{\ssf}{\textsf}
\newcommand{\mb}{}

\newcommand{\SP}{\text{Spec }}

\newcommand{\Kbm}[2]{\overline{\mathcal{M}}_{#1}(#2)}
\newcommand{\kbm}[2]{\overline{M}_{#1}(#2)}

\newcommand{\marpar}[1]{}

\newcommand{\lt}{\left}
\newcommand{\rt}{\right}

\newcommand{\mc}{\mathcal}

\newcommand{\OO}{\mathcal O}

\def\PP{{\mathbb P}}

\begin{document}

\title[Rationally simply-connected hypersurfaces]
{Hypersurfaces of low degree are rationally simply-connected}
\
\author[Starr]{Jason Michael Starr}
\address{Department of Mathematics \\
  Massachusetts Institute of Technology \\ Cambridge MA 02139}
\email{jstarr@math.mit.edu} 

\date{\today}

\begin{abstract}
For a general hypersurface of degree $d$ in projective $n$-space, if
$n \geq d^2$ the spaces of $2$-pointed rational curves on the
hypersurface are rationally connected; thus the hypersurfaces are
\emph{rationally simply connected}.
This paper proves stronger versions of theorems in ~\cite{HS2}.
\end{abstract}

\maketitle


\section{Introduction} \label{sec-intro}
\marpar{sec-intro}

\noindent
A smooth projective variety is \emph{rationally connected} if every
pair of points is contained in a rational curve, cf. ~\cite{K},
~\cite{De}.  This is analogous to 
the notion of path-connectedness in topology.  A path-connected
topological space is simply-connected if the space of based paths is
path-connected.  By analogy, Mazur suggested definining
a rationally connected variety to be
\emph{rationally simply-connected} if the space of based,
2-pointed rational curves of fixed homology class is rationally
connected.  Unfortunately, this notion is a bit too strong.  The
condition should be imposed only if the homology class is suitably positive
and if the basepoints are suitably general.  Also, since the space of
2-pointed rational curves is typically not compact, it is compactified
using the Kontsevich moduli space, cf. ~\cite{FP}.

\begin{thm} \label{thm-main}
In characteristic $0$, a general hypersurface of degree $d$ in $\PP^n$
is rationally simply connected if $n\geq d^2$ and $d\geq 2$.  
Precisely, for every $e\geq 2$, the evaluation morphism
$\text{ev}:\Kbm{0,2}{X,e} \rightarrow X \times X$ is surjective and a
general fiber is irreducible, reduced and rationally connected.  This
also holds for $(n,d) = (\geq 2,1)$.  
\end{thm}

\medskip\noindent
Tsen's theorem proves that a Fano hypersurface defined over the
function field of a curve has a rational point, cf. ~\cite{Tsen36}. 
Rational connectedness, introduced by Campana ~\cite{Ca} and
Koll\'ar-Miyaoka-Mori ~\cite{KMM}, is the natural notion for
generalizing Tsen's theorem, cf. ~\cite{GHS}, ~\cite{dJS} and
~\cite{GHMS}.  The analogue of Tsen's theorem over a higher-dimensional
base is the Tsen-Lang theorem, cf. ~\cite{Lang52}.  To generalize the
Tsen-Lang theorem, the correct notion should be some version of 
\emph{higher rational connectedness}

\medskip\noindent
In truth, the best definition of rational
simple-connectedness is not yet clear.  At issue is a precise
definition of sufficiently positive homology class.
However, if
$\text{Pic}(X) = \ZZ$ there is only one possible definition: a
collection of homology classes $\beta \in H_2(X,\ZZ)$ is
\emph{sufficiently positive} if for one ample invertible sheaf $L$,
and hence every ample invertible sheaf, there exists an integer
$e_0=e_0(L)$ such that $\beta$ is in the collection if $c_1(L)\cdot
\beta \geq e_0$.  This is the notion that appears in
Theorem~\ref{thm-main}.  

\begin{prin}[de Jong] \label{prin-deJong}
\marpar{prin-deJong}
A proper, smooth variety defined over the function field of a surface
over an algebraically closed field of characteristic $0$ has a
rational point if the base-change to the algebraic closure of the
function field is rationally simply-connected (in a slightly stronger
sense than used here)  
and if a certain Brauer
obstruction vanishes.
\end{prin}

\medskip\noindent
In proving the principle, there are
technical difficulties related to singularities of
spaces of rational curves, the geometric analogue of weak
approximation, and \emph{the} Brauer obstruction (as opposed to
\emph{a} Brauer obstruction).
In a forthcoming paper with
A. J. de Jong, de Jong's strategy for proving
Principle~\ref{prin-deJong} is given, and Principle~\ref{prin-deJong}
is proved under some additional hypotheses. 

\medskip\noindent
\textbf{Acknowledgments:}  A re-investigation of ~\cite{HS2} leading
to this paper was motivated by Theorem \ref{thm-uniruled} of A. J. de
Jong.  I am grateful to de Jong for many useful conversations.

\section{Sketch of the proof} \label{sec-sketch}
\marpar{sec-sketch}

\noindent
The proof of Theorem~\ref{thm-main} uses ideas from ~\cite{HS2}.
It is an induction argument.  The induction step
uses the fiber-by-fiber connected sum of 2 families of
marked rational curves.  
Three new results go into the proof.  First, very
  twisting morphisms are systematically studied leading to a simpler
  version of the induction argument, Proposition ~\ref{prop-induction}. 
The other two results concern the base case for the induction.
Proposition~\ref{prop-conics} proves Theorem~\ref{thm-main} in
case $e=2$.  The last result 
is existence of a very twisting morphism 
$\zeta:\PP^1 \rightarrow \Kbm{0,1}{X,1}$.  Some notation is necessary
to state the result.

\begin{notat} \label{notat-xh}
\marpar{notat-xh}
Let $n\geq 2$.  The scheme $\Kbm{0,1}{\PP^n,1}$ is the partial flag
variety $\text{Flag}(1,2;n+1)$.  There are morphisms
$\text{ev}:\Kbm{0,1}{\PP^n,1} \rightarrow \PP^n$ and
$\text{pr}:\Kbm{0,1}{\PP^n,1} \rightarrow \Kbm{0,0}{\PP^n,1} =
\text{Grass}(2,n+1)$.  

\noindent
Denote by $h$ the divisor class $\text{ev}^*
C_1(\OO_{\PP^n}(1))$.  Denote by $x$ the divisor class $\text{pr}^*
C_1(\OO_{\text{Grass}}(1))$, where $\OO_{\text{Grass}}(1)$ is the
invertible sheaf giving the Pl\"ucker embedding.  

\noindent
The NEF cone of
$\Kbm{0,1}{\PP^n,1}$ is the set of divisor classes $\ZZ_{\geq
  0}\{x,h\}$.  
\end{notat}

\medskip\noindent
The degree of a very twisting
morphism $\zeta$
is large:  
If $n=d^2$ then $\text{deg}(\zeta^* h)
\geq d^2-d-1$ and $\text{deg}(\zeta^* x) \geq 2(d^2-d-1)$.
So it is unreasonable to try to directly construct a very twisting
morphism.  
The strategy is to instead consider a nodal, reducible, genus $0$ curve $C$,
and a morphism
$\zeta:C\rightarrow \Kbm{0,1}{\mb{X},1}$ whose restriction to 
every irreducible component of $C$ is a
twisting morphism of minimal degree. 

\medskip\noindent
The curve $C$ deforms to a smooth genus $0$ curve, i.e., to $\PP^1$.
If $C$ has many irreducible components, 
a deformation of $\zeta$ to a morphism from $\PP^1$ may be a very
twisting morphism.
Proposition~\ref{prop-q}
gives a criterion for this.
The essential case is $n=d^2$.  Then the criterion is
surjectivity of the derivative map between the Zariski tangent spaces
of two parameter spaces.  Sections \ref{sec-deriv}, \ref{sec-proof},
\ref{sec-pfmain}, 
and particularly 
Proposition~\ref{prop-twistlines}, prove the
derivative map is surjective.  This 
is the heart of the article and consists
of an involved deformation theory analysis.

\section{Twisting and very twisting morphisms}
\label{sec-twist}
\marpar{sec-twist}

\begin{hyp} \label{hyp-K}
\marpar{hyp-K}
Let $K$ be an algebraically closed field.  All schemes, algebraic
spaces and Deligne-Mumford stacks are defined over
$K$, and all morphisms and $1$-morphisms commute with the morphisms to
$\SP K$.
\end{hyp}

\medskip\noindent
Every rank $r$, locally free $\OO_{\PP^1}$-module is isomorphic to a
direct sum of invertible sheaves by Grothendieck's lemma
~\cite[Exer. V.2.6]{H}, i.e., 
$\mc{E} \cong \OO_{\PP^1}(a_1)\oplus
\dots \oplus \OO_{\PP^1}(a_r)$ for a sequence of integers $a_1\leq
\dots \leq a_r$.    

\begin{defn} \label{defn-null}
\marpar{defn-null}
The \emph{negativity of $\mc{E}$} is
the largest integer $n$ (possibly $0$) such that
$a_i<0$ for all $i\leq n$.  The \emph{nullity of $\mc{E}$} is
the largest integer
$z$ (possibly $0$) such that $a_i=0$ for all $n<i\leq n+z$.  The
\emph{positivity of $\mc{E}$} is  
the difference $p:= r-n-z$.  

\noindent
If $n(\mc{E})=0$, then $\mc{E}$ is
\emph{generated by global sections}.  If $n(\mc{E})+z(\mc{E})=0$ and
$p(\mc{E}) > 0$, then $\mc{E}$ is \emph{ample}.
\end{defn}

\medskip\noindent
Let $\pi:Y\rightarrow Z$
be a morphism of finite type Deligne-Mumford stacks.  Denote by
$Z^0$ the maximal open substack of $Z$ that is smooth over $K$.
Denote by $Y^0$ the maximal open substack of $\pi^{-1}(Z^0)$ on which
$\pi$ is smooth.  Denote by $T_\pi$ the locally free sheaf on $Y^0$
that is the dual of the sheaf of relative differentials of $\pi$,
i.e., $T_\pi$ is the vertical tangent bundle of $\pi$. 

\begin{defn} \label{defn-relfree}
\marpar{defn-relfree}
A non-constant 
morphism $f:\PP^1 \rightarrow Y$ is \emph{$\pi$-relatively free},
resp. \emph{$\pi$-relatively very free}, if
\begin{enumerate}
\item[(i)]
$f(\PP^1) \subset Y^0$,
\item[(ii)]
$f^*\pi^* T_{Z^0}$ is generated by global sections, and
\item[(iii)]
$f^* T_{\pi}$ is generated by global sections (resp. ample).
\end{enumerate}
For a $\pi$-relatively free morphism $f$, the \emph{nullity of $f$} is
$z(f^* T_\pi)$, and the \emph{positivity of $f$} is 
$p(f^* T_\pi)$.
\end{defn}

\begin{lem} \label{lem-relfree}
\marpar{lem-relfree}
The smooth locus of the 1-morphism
$\text{Hom}(\PP^1,\pi): \text{Hom}(\PP^1,Y)\rightarrow
\text{Hom}(\PP^1,Z)$ contains the locus of $\pi$-relatively free morphisms.
\end{lem}

\begin{proof}
Associated to a $\pi$-relatively free morphism $f$
there is a short exact sequence, 
$$
\begin{CD}
0 @>>> f^* T_\pi @>>> f^* T_Y @> d\pi >> f^*\pi^* T_Z @>>> 0.
\end{CD}
$$
Because the outer terms are generated by global sections, in
particular $h^1(\PP^1,f^* T_\pi)=h^1(\PP^1,f^* T_Y) = h^1(\PP^1, f^*
\pi^* T_Z) = 0$.  Therefore both $[f]\in \text{Hom}(\PP^1,Y)$ and
$[\pi\circ f]\in \text{Hom}(\PP^1,Z)$ are smooth points, 
and $H^0(\PP^1,d\pi):H^0(\PP^1,f^* T_Y) \rightarrow H^0(\PP^1,f^*
\pi^* T_Z)$ is surjective.  Because this is the derivative of
$\text{Hom}(\PP^1,\pi)$ at $[f]$, the Jacobian criterion implies
$\text{Hom}(\PP^1,\pi)$ is smooth at $[f]$.
\end{proof}

\begin{rmk} \label{rmk-relfree}
\marpar{rmk-relfree}
The smooth locus of
$\PP^1 \times 
\text{Hom}(\PP^1,Y) \rightarrow Y$ contains $\PP^1\times \{f\}$ for
every $\pi$-relatively free morphism $f$.  For every open subset
$U\subset Y$ whose closure intersects $\text{Image}(f)$, every
sufficiently small
deformation $f_\epsilon$ of $f$ is $\pi$-relatively free and
$\text{Image}(f_\epsilon)$ intersects $U$.
\end{rmk}

\medskip\noindent
The following proposition is the relative version of the criterion
~\cite[Theorem IV.5.8]{K}.  

\begin{prop} \label{prop-relfree}
\marpar{prop-relfree}
Let $\pi:Y \rightarrow Z$ be a proper morphism of irreducible, finite type
Deligne-Mumford stacks, and denote by $|\pi|$ the induced morphism of
coarse moduli spaces.
\begin{enumerate}
\item[(i)]
If there exists a $\pi$-relatively free morphism $f$ with positivity
$>0$, 
then every irreducible component of the geometric generic
fiber of $|\pi|$ is uniruled.
\item[(ii)]
If also $\text{char}(K)=0$, the dimension of the MRC
quotient of the geometric generic fiber of $|\pi|$ 
is at most the nullity of $f$,
i.e., the fiber dimension of the MRC quotient morphism is at least the
positivity of $f$. 
\end{enumerate}
\end{prop}

\begin{proof}
\textbf{(i):}
This is essentially ~\cite[Lem. 4.4]{GHMS}, which in turn relies on
the \emph{Rigidity Lemma}, cf. ~\cite[p. 43]{AVar}. 
Denote $Y_f = \PP^1 \times_{\pi\circ f,Z,\pi} Y$
and denote the projection by 
$\pi_f:Y_f \rightarrow \PP^1$.  The morphism $f$ uniquely determines 
a section $\sigma$ of $\pi_f$.  The projection
$\pi_f$ is smooth along $\sigma(\PP^1)$.  The sheaf $\sigma^*
T_{\pi_f} = f^* T_\pi$ is generated by global sections.

\medskip\noindent
Denote by $M \subset \text{Hom}(\PP^1,Y_f)$ the locally closed
substack of the Hom stack parametrizing morphisms
$\tau:\PP^1\rightarrow Y_f$ with $\pi_f\circ \tau =
\text{Id}_f$.
By
~\cite[Prop. II.3.5]{K}, the smooth locus of 
$\text{ev}: \PP^1 \times M \rightarrow Y_f$ contains $\PP^1
\times \{[\sigma]\}$, and the fiber dimension of $\text{ev}$ is
$h^0(\PP^1,\sigma^* T_{\pi_f}(-1))$.  The proof is unchanged in the
case of Deligne-Mumford stacks.  

\medskip\noindent
Because the positivity of $f$ is
$>0$, the fiber dimension of $\text{ev}$ is $>0$.  
By ~\cite[Lem. 4.4]{GHMS}, for
every $t\in \PP^1$ there exists a rational curve containing
$\sigma(t)$ and contained in 
$|\pi_f|^{-1}(t)$.
The hypothesis $K=\CC$ from ~\cite{GHMS} 
is not used in the proof of ~\cite[Lem. 4.4]{GHMS}.  

\medskip\noindent
For every $t\in \PP^1$ there is a rational curve containing $f(t)$
contained in a
fiber of $|\pi|$. 
By Remark~\ref{rmk-relfree}, the morphism
$\PP^1 \times \text{Hom}(\PP^1,Y) \rightarrow Y$ is smooth at
$(t,[f])$. So there is a dense open subset $U$ of $Y^0$ such that for
every geometric point $y\in U$, there is a rational curve containing
$y$ and contained in a fiber
of $|\pi|$. 

\medskip\noindent
\textbf{(ii):}
By ~\cite{Hir}, assume $Y$ and $Z$ are smooth and irreducible.  
The relative MRC quotient for Deligne-Mumford stacks is a datum
$(Q,\rho,U,\psi)$ of 
a smooth, irreducible algebraic space $Q$, a proper morphism
$\rho:Q\rightarrow Z$ whose geometric generic fiber is not uniruled,
an open subset $U\subset Y$ whose complement has codimension $\geq
2$, and a generically smooth 1-morphism $\psi:U \rightarrow Q$ such that
$\pi|_U$ is 2-equivalent to $\rho\circ \psi$.  Existence of the
relative MRC quotient for Deligne-Mumford stacks is proved in
~\cite{Sfree}.  

\medskip\noindent
Because $f:\PP^1 \rightarrow Y$ is $\pi$-relatively free, it is free.
By ~\cite[Prop. II.3.7]{K}, the morphism $f$ can be deformed so that
$f(\PP^1)$ is contained in $U\cap Y^0$ and $f(\PP^1)$ intersects the smooth
locus of $\psi$.   
By the Jacobian criterion, at every point $y\in Y^0$ the derivative
map $d\pi_y: T_{Y}\otimes \kappa(y) \rightarrow T_{Z}\otimes
\kappa(\pi(y))$ is surjective.  The derivative map
$d\rho_{\psi(y)}: T_{Q}\otimes \kappa(\psi(y)) \rightarrow T_Z \otimes
\kappa(\pi(y))$ factors $d\pi_y$.  So $d\rho_{\psi(y)}$ is also
surjective.  By the Jacobian criterion, $\rho$ is smooth at $\psi(y)$,
i.e., $\rho$ is smooth along $\psi(Y^0)$.  In particular, $\rho$ is
smooth along $\psi(f(\PP^1))$.

\medskip\noindent
There is a sheaf homomorphism $T_\pi|_U
\rightarrow \psi^* T_\rho$ whose restriction to the smooth locus of
$\psi$ is surjective.
So there is a sheaf homomorphism $\alpha:f^* T_\pi \rightarrow f^*\psi^*
T_{\rho}$ whose cokernel is a torsion sheaf.
Because $f^* T_{\pi}$ is generated by
global sections, also $\text{Image}(\alpha)$ and $f^* \psi^* T_{\rho}$
are generated by global sections.
Since the geometric generic fiber of $\rho$ is not uniruled, by
(i) the positivity of $\psi\circ f$ is $0$.  This implies that
$\text{Image}(\alpha) = f^*\psi^* T_{\rho}$ and that
$f^*\psi^*
T_\rho$ is a quotient of $f^* T_\pi/H^0(\PP^1,f^*T_\pi(-1))\otimes
\OO_{\PP^1}(1)$.  Therefore the rank of $f^* \psi^* T_{\rho}$
is at most the nullity of $f$.
In other words, the dimension of the geometric generic fiber of
$\rho$ is at most the nullity of $f$.  
\end{proof}

\medskip\noindent
The main definition of both ~\cite{HS2} and this paper is the
following.  It is slightly different than in ~\cite{HS2} because it is
used differently here.

\begin{defn} \label{defn-twist}
\marpar{defn-twist}
Let $X$ be a
quasi-projective scheme, and denote by $X^0$ the smooth
locus of $X$.  Let $r\geq 0$ and $e >0$ be integers.
Let $\zeta:\PP^1 \rightarrow \Kbm{0,r}{X^0,e}$ be a 1-morphism.
The 1-morphism $\zeta$ is \emph{twisting},
resp. \emph{very twisting}, if
\begin{enumerate}
\item[(i)]
the morphism $\text{ev}:\Kbm{0,r}{\mb{X}^0,e} \rightarrow
(\mb{X}^0)^r$ is unobstructed at every geometric point of the image of
$\zeta$, 
\item[(ii)]
the morphism $\zeta$ is $\text{ev}$-relatively free,
resp. $\text{ev}$-relatively very free, 
\item[(iii)]
for every $i=1,\dots,r$, the degree of $\zeta^* \psi_i$ is
nonpositive, and
\item[(iv)]
the image under $\zeta$ of the geometric generic point of $\PP^1$ is a
stable map with irreducible domain.
\end{enumerate}
Let $T$ be a Deligne-Mumford stack and let $\zeta:\PP^1
\times T \rightarrow \Kbm{0,r}{\mb{X}^0,e}$ be a 1-morphism.  The
morphism $\zeta$ is \emph{twisting relative to $T$}, resp. \emph{very
  twisting relative to $T$}, if for every geometric point of $T$ the
restriction of $\zeta$ over this point is twisting, resp. very twisting.
\end{defn}

\begin{rmk} \label{rmk-twist}
\marpar{rmk-twist}
\begin{enumerate}
\item[(i)]
If $r=0$, the morphism $\text{ev}$ is just the structure morphism to
$\SP(K)$.  So in this case, $\zeta$ is $\text{ev}$-relatively free,
resp. $\text{ev}$-relatively very free, iff it is free, resp. very free.
\item[(ii)]
Condition (i) implies that $\text{ev}$ is smooth along the image of
$\zeta$.  If $\text{char}(K)=0$ and $\mb{X}$ is  
a general hypersurface of degree $d< \frac{n+1}{2}$ in $\PP^n$, the
smooth locus of $\text{ev}$ is precisely the set of points where
$\text{ev}$ is unobstructed.  But for some schemes 
the unobstructed
locus is strictly smaller than the smooth locus, e.g., for a
general hypersurface of degree $n$ in $\PP^n$ and $e=4$.
\item[(iii)]
The 1-morphism $\zeta$ is equivalent to a datum
$((\pi:\Sigma\rightarrow \PP^1,\sigma_1,\dots,\sigma_r),g:\Sigma
\rightarrow \mb{X}^0)$.  The class $\zeta^* \psi_i$ is simply the
divisor class of $\sigma_i^* \OO_{\Sigma}(-\sigma_i(\PP^1))$.
Therefore condition (iii) states that the self-intersection 
$(\sigma_i(\PP^1)\cdot \sigma_i(\PP^1))_\Sigma$ is nonnegative for every
$i=1,\dots,r$.
\item[(iv)]
Often the morphism $(\pi,g):\Sigma \rightarrow \PP^1\times \mb{X}^0$
is unramified and is \'etale locally a regular embedding, i.e., the
sheaf homomorphism $d(\pi,g)^\dagger:(\pi,g)^* \Omega_{\PP^1\times \mb{X}^0}
\rightarrow \Omega_\Sigma$ is surjective and the kernel is locally
free.  In particular, this is true if $e=1$. If $d(\pi,g)^\dagger$ is
surjective and the kernel is locally free, 
denote by $N_{(\pi,g)}$ the dual of the kernel.  The
morphism $\text{ev}$ is unobstructed along $\zeta(\PP^1)$ iff
$R^1\pi_*(N_{(\pi,g)}(-(\sigma_1(\PP^1)+\dots+\sigma_r(\PP^1))))$ is
$0$ and then $\zeta^* T_{\text{ev}}$ is the locally free
sheaf $\pi_*(N_{(\pi,g)}(-(\sigma_1(\PP^1)+\dots+\sigma_r(\PP^1))))$. 
\item[(v)]
By Remark~\ref{rmk-relfree}, if there exists a 1-morphism $\zeta:\PP^1
\rightarrow \Kbm{0,r}{\mb{X},e}$ that satisfies (i)--(iii), and if
$\zeta(\PP^1)$ intersects the closure of the open set parametrizing
stable maps with irreducible domain, then a small deformation of
$\zeta$ is a 1-morphism satisfying (i)--(iv).
\end{enumerate}
\end{rmk}

\medskip\noindent
If $r\geq 2$ consideration of the Picard group of $\Sigma$
implies that $\text{deg}(\zeta^* \psi_i) = 0$ for all $i=1,\dots,r$.
However if $r=1$ -- the main case of interest -- then $-\text{deg}(\zeta^*
\psi)$ can be arbitrarily positive.

\begin{prop} \label{prop-psineg}
\marpar{prop-psineg}
Let $\zeta_0:\PP^1 \rightarrow \Kbm{0,1}{\mb{X}^0,e}$ be a very
twisting morphism.  Denote $a_0 = -\text{deg}(\zeta_0^*\psi)$.  
\begin{enumerate}
\item[(i)]
If $a_0$ is odd, there exists an integer $a_1$ and
for every integer $a\geq a_1$ there is a very twisting morphism
$\zeta_a:\PP^1  
\rightarrow \Kbm{0,1}{\mb{X}^0,e}$ such that $a=
-\text{deg}(\zeta_a^*\psi)$.  
\item[(ii)]
If $a_0$ is even, there exists an integer $a_1$ and for
every even integer $a\geq a_1$ there is a very twisting morphism
$\zeta_a:\PP^1 \rightarrow \Kbm{0,1}{\mb{X}^0,e}$ such that $a=
-\text{deg}(\zeta_a^* \psi)$.
\end{enumerate}
\end{prop} 

\begin{proof}
\textbf{(i):}  Since $\zeta_0$ is very twisting, $\zeta_0^*
T_{\text{ev}} \cong \OO_{\PP^1}(b_1) + \dots + \OO_{\PP^1}(b_t)$ for
integers $1\leq b_1 \leq \dots \leq b_t$.  Define,
$$
a_1 = 2a_0 \lt\lceil \frac{a_0+b_1}{2b_1} \rt\rceil.  
$$
Let $a \geq a_1$ be an integer.  There are
2 cases depending on whether $a$ is even or odd.

\medskip\noindent
If $a$ is even, then $a=q(2a_0) + r$ for an integer $q$ and an integer
$r$ satisfying $0\leq r < 2a_0$.  Because $a$ is even, $r=2r'$ for an
integer $r'$ satisfying $0\leq r' < a_0$.  Define $m=2q$.  
Because $a \geq a_1$, $m
\geq \frac{a_0+b_1}{b_1}$, and thus $mb_1 - r' > b_1
>0$.  

\medskip\noindent
If $a$ is odd, then $a+a_0 = q(2a_0) + r$ for an integer $q$ and an
integer $r$ satisfying $0\leq r < 2a_0$.  Because $a+a_0$ is even,
$r=2r'$ for an integer $r'$ satisfying $0 \leq r' < a_0$.  Define
$m=2q-1$.  Because $a\geq a_1$, $m \geq \frac{a_0}{b_1}$, and thus
$mb_1 - r' > 0$.

\medskip\noindent
In each case, $a=ma_0 + 2r'$ where $0\leq r' < a_0$ and where $mb_1 -
r' > 0$.  Let $h:\PP^1 \rightarrow \PP^1$ be a finite morphism of
degree $m$.  The morphism $\zeta_0\circ h:\PP^1 \rightarrow
\Kbm{0,1}{\mb{X},e}$ is very twisting.  Moreover, $(\zeta_0\circ h)^*
T_{\text{ev}} \cong \OO_{\PP^1}(mb_1) \oplus \dots \oplus
\OO_{\PP^1}(mb_t)$ and $-\text{deg}((\zeta_0\circ h)^* \psi) = ma_0$.

\medskip\noindent
Denote by $((\pi:\Sigma\rightarrow \PP^1, \sigma:\PP^1 \rightarrow
\Sigma),g:\Sigma \rightarrow \mb{X})$ the datum giving rise to
$\zeta_0\circ h$.  There exists a section $\sigma':\PP^1 \rightarrow
\Sigma$ such that the divisor $\sigma'(\PP^1) \subset \Sigma$ is in
the linear equivalence class of $|\sigma(\PP^1) + r' \pi^{-1}(0)|$.
Denote by $\zeta_a:\PP^1 \rightarrow \Kbm{0,1}{\mb{X},e}$ the
1-morphism associated to the datum $((\pi:\Sigma\rightarrow \PP^1,
\sigma':\PP^1 \rightarrow \Sigma),g:\Sigma \rightarrow \mb{X})$.  

\medskip\noindent
For every $\ssf{s}\in \PP^1$, the obstruction group of $\text{ev}$ at
$\zeta_0(h(\ssf{s}))$ is the hypercohomology group
$$
\mathbf{H}^2(\pi^{-1}(\ssf{s}), L^\vee_{(\pi,g)} \otimes
\OO_{\pi^{-1}(\ssf{s})}(-\sigma(\ssf{s}))),
$$
where $L_{(\pi,g)}$ is the cotangent complex of the morphism
$(\pi,g):\Sigma \rightarrow \PP^1 \times \mb{X}$ and where
$L^\vee_{(\pi,g)}$ is the object
$R\textit{Hom}_{\OO_\Sigma}(L_{(\pi,g)},\OO_\Sigma)$ in the derived
category of quasi-coherent sheaves on $\Sigma$.  
By hypothesis, the obstruction group is $0$.
By construction, $\OO_{\pi^{-1}(\ssf{s})}(-\sigma'(\ssf{s})) \cong
\OO_{\pi^{-1}(\ssf{s})}(-\sigma'(\ssf{s}))$.  Therefore, the
obstruction group of $\text{ev}$ at $\zeta_a(\ssf{s})$ is $0$.
So $\zeta_a$ satisfies (i) of Definition~\ref{defn-twist}.

\medskip\noindent
Moreover, $\zeta_a^* T_{\text{ev}}$ is the cohomology sheaf,
$$
\mathcal{H}^1 R\pi_*(
L^\vee_{(\pi,g)}(-\sigma'(\PP^1)) ) = \mathcal{H}^1 R\pi_*(
L^\vee_{(\pi,g)}(-\sigma(\PP^1)) )\otimes \OO_{\PP^1}(-r'),
$$
i.e., $h^*\zeta_0^* T_{\text{ev}} \otimes \OO_{\PP^1}(-r')$.
So $\zeta_a^* T_{\text{ev}} \cong \OO_{\PP^1}(mb_1 - r') \oplus \dots
\oplus \OO_{\PP^1}(mb_t - r')$.  By construction, $mb_1 - r' >0$.
Therefore $\zeta_a^* T_{\text{ev}}$ is ample.  So $\zeta_a$ satisfies
(ii) of Definition~\ref{defn-twist} for a very twisting morphism.

\medskip\noindent
Finally, $(\sigma'(\PP^1)\cdot \sigma'(\PP^1))_\Sigma =
(\sigma(\PP^1)\cdot \sigma(\PP^1))_\Sigma + 2r'$.  Therefore,
$$
-\text{deg}(\zeta_a^* \psi) = -\text{deg}(g^* \zeta_0^* \psi) + 2r' =
ma_0 + 2r' = a.
$$
So $\zeta_a$ satisfies (iii) of Definition~\ref{defn-twist} for a 
very twisting morphism.  Because $\zeta_0$ satisfies $(iv)$ of
Definition~\ref{defn-twist}, also $\zeta_a$ satisfies $(iv)$.
Therefore $\zeta_a$ is very twisting and
$-\text{deg}(\zeta_a^*\psi) = a$.

\medskip\noindent
\textbf{(ii):}
This is exactly as in the even case of (i).  The odd case of (i) does
not work because $a+a_0$ is not even.
\end{proof}

\medskip\noindent
The main result about very twisting morphisms is the following
theorem, which forms the induction step in the proof of
Theorem~\ref{thm-main}.

\begin{prop} \label{prop-induction}
\marpar{prop-induction}
Let $\mb{X}$ be a quasi-projective variety.  Let $e_1,e_2 > 0$ be
integers.  Let $r$ equal $1$, resp. $2$.  Let $\zeta_1:\PP^1
\rightarrow \Kbm{0,r}{\mb{X},e_1}$ be an $\text{ev}$-relatively very
free morphism mapping a general point of $\PP^1$ to a stable map with
irreducible domain and with $\text{ev} \circ \zeta_1$ nonconstant
(resp. $\text{ev}_1\circ \zeta_2$ nonconstant). 
Let $\zeta_2:\PP^1 \rightarrow \Kbm{0,1}{\mb{X},e_2}$ be a very
twisting morphism.

\medskip\noindent  
Assume that $\Kbm{0,0}{\mb{X},\epsilon}$ is irreducible for all positive
integers $\epsilon$ that are sufficiently divisible.  Assume one of
the following,
\begin{enumerate}
\item[(i)]
$-\text{deg}(\zeta^*_1\psi) > 0$
(resp. $-\text{deg}(\zeta^*_1\psi_1)>0$), or 
\item[(ii)]
for every general
degree $e_1$ morphism $g:\PP^1 \rightarrow \mb{X}$ there exists a
twisting morphism $\zeta_g:\PP^1 \rightarrow \Kbm{0,1}{\mb{X},e_2}$
such that $\text{ev}\circ \zeta_g = g$.  
\end{enumerate}
Then there exists an
$\text{ev}$-relatively very free morphism $\zeta:\PP^1 \rightarrow
\Kbm{0,r}{\mb{X},e_1+e_2}$ mapping a general point of $\PP^1$ to a
stable map with irreducible domain and 
with $\text{ev}\circ \zeta_1$ nonconstant
(resp. $\text{ev}_1\circ \zeta_2$ nonconstant). 

\medskip\noindent
Moreover, if
$r=1$ and $\zeta_1$ is very twisting, then there exists such a $\zeta$
that is very twisting and has $-\text{deg}(\zeta^*\psi)=0$.  
\end{prop}

\begin{proof}
To ease notation, if $r=1$ denote $\text{ev}$ also by $\text{ev}_1$.
By Proposition~\ref{prop-psineg}, assume that
$-\text{deg}(\zeta_2^*\psi) > 0$.  If $r=1$ and $\zeta_1$ is very
twisting, also assume that $-\text{deg}(\zeta_1^*\psi) > 0$, i.e., (i)
applies.  

\medskip\noindent
Let $\epsilon$ be a positive integer such
that $\Kbm{0,0}{\mb{X},\epsilon}$ is irreducible.  
After precomposing $\zeta_1$ and $\zeta_2$ with finite morphism $\PP^1
\rightarrow 
\PP^1$, assume that
$\text{deg}(\text{ev} \circ \zeta_1) = \text{deg}(\text{ev} \circ
\zeta_2)=\epsilon$.
Denote by $\text{Hom}(\PP^1,\mb{X})_\epsilon$ the open subset of
$\text{Hom}(\PP^1,\mb{X})$ parametrizing morphisms of degree
$\epsilon$; this is irreducible by hypothesis.  There are morphisms,
$$
\begin{CD}
\text{Hom}(\PP^1,\text{ev}_1): &
\text{Hom}(\PP^1,\Kbm{0,r}{\mb{X},e_1}) @>>>
\text{Hom}(\PP^1,\mb{X}) \\
\text{Hom}(\PP^1,\text{ev}): & \text{Hom}(\PP^1,\Kbm{0,1}{\mb{X},e_2})
@>>> \text{Hom}(\PP^1,\mb{X})
\end{CD}
$$
By Lemma~\ref{lem-relfree}, $\text{Hom}(\PP^1,\text{ev}_1)$ is smooth
at $[\zeta_1]$ and $\text{Hom}(\PP^1,\text{ev})$ is smooth at
$[\zeta_2]$.  So the image of each is a dense open subset of
$\text{Hom}(\PP^1,\mb{X})_\epsilon$.  Therefore, after deforming
$\zeta_1$ and $\zeta_2$, assume that also $\text{ev}_1\circ \zeta_1 =
\text{ev}\circ \zeta_2$.

\medskip\noindent
Denote by $((\pi_2:\Sigma_2 \rightarrow \PP^1, \sigma_{2,1}),g_2)$ the
datum associated to $\zeta_2$.
If $r=1$, denote by $((\pi_1:\Sigma_1\rightarrow
\PP^1,\sigma_{1,1}),g_1)$ the datum associated to $\zeta_1$.  If
$r=2$, denote by $((\pi_1:\Sigma_1 \rightarrow
\PP^1,\sigma_{1,1},\sigma_{1,2}),g_1)$ the datum associated to
$\zeta_1$. 
Denote the two cases in the statement as \emph{Case (i)} and
\emph{Case (ii)} respectively.
In Case (ii) a similar argument to the last paragaph proves, after
deforming $\zeta_1$ and $\zeta_2$ further, there exists a dense open
subset $U\subset \PP^1$ such that for every closed point $\ssf{s}\in
U$, 
\begin{enumerate}
\item[(i)]
$\pi_1^{-1}(\ssf{s})$ is irreducible, and 
\item[(ii)]
there exists a twisting
morphism $\zeta_{\ssf{s}}:\pi_1^{-1}(\ssf{s}) \rightarrow
\Kbm{0,1}{\mb{X},e_2}$ such that
$\zeta_{\ssf{s}}(\sigma_{1,1}(\ssf{s})) = \zeta_2(\ssf{s}).$
\end{enumerate}

\medskip\noindent
In Case (i), denote $C=C_0=\sigma_{1,1}(\PP^1)$.  In Case (ii), denote
$C_0=\sigma_{1,1}(\PP^1)$.  Let $\delta > 0$ be an integer such that
$(C_0\cdot C_0)_{\Sigma_1} + \delta > 0$.  Let
$\ssf{s}_1,\dots,\ssf{s}_\delta \in U$ be distinct points.  For
$i=1,\dots,\delta$, denote $C_i = \pi_1^{-1}(\ssf{s}_i)$.  Denote
$C=C_0 \cup C_1 \cup \dots \cup C_\delta$.  In each case, the linear
system $|C|$ on $\Sigma_1$ contains an irreducible curve that
intersects $\sigma_{1,1}(\PP^1)$ (resp. $\sigma_{1,1}(\PP^1)\cup
\sigma_{1,2}(\PP^1)$) transversely in finitely many points.  Let
$\mc{D} \subset \PP^1_t \times \Sigma_1$ be a divisor such that
$\text{pr}_{\PP^1}^{-1}(0)\cap \mc{D} \subset \Sigma_1$ equals $C$ and
such that for general $\ssf{t}\in \PP^1_t$,
$\text{pr}_{\PP^1}^{-1}(\ssf{t}) \cap \mc{D} \subset \Sigma_1$ is an
irreducible divisor that intersects $\sigma_{1,1}(\PP^1)$
(resp. $\sigma_{1,1}(\PP^1)\cup \sigma_{1,2}(\PP^1)$) transversely in
finitely many points (the subscript ``t'' in $\PP^1_t$ is to
distinguish $\PP^1_t$ from the target of $\pi_1$).  

\medskip\noindent
Denote by $\zeta_C:C\rightarrow \Kbm{0,1}{\mb{X},e_2}$ the unique
morphism such that $\zeta_C\circ \sigma_{1,1} = \zeta_2$ and such that
$\zeta_C|_{C_i} = \zeta_{\ssf{s}_i}$ for $i=1,\dots,\delta$ in Case
(ii).  There is a 1-morphism of relative Hom stacks,
$$
\widehat{\text{ev}}: \text{Hom}_{\PP^1_t}(\mc{D},
\PP^1_t \times \Kbm{0,1}{\mb{X},e_2}) \rightarrow
\text{Hom}_{\PP^1_t}(\mc{D}, \PP^1_t \times \mb{X}).
$$
There is a point $(0,[\zeta_C]) \in
\text{Hom}_{\PP^1_t}(\mc{D},\PP^1_t \times \Kbm{0,1}{\mb{X},e_2})$
lying over $0\in \PP^1_t$.  By the same argument as in the proof of
Lemma~\ref{lem-relfree}, the morphism $\widehat{\text{ev}}$ is smooth
at $(0,[\zeta_C])$.  

\medskip\noindent
The following morphism defines a section $\widehat{g}_1$ of
$\text{Hom}_{\PP^1_t}(\mc{D},\PP^1_t \times \mb{X}) \rightarrow \PP^1_t$,
$$
(\text{pr}_{\PP^1_t},
g_1\circ \text{pr}_{\Sigma_1}):\mc{D} \rightarrow \PP^1_t\times \mb{X}.
$$
Denote by $\text{pr}_{\PP^1_t}:H \rightarrow \PP^1_t$ the fiber
product,
$$
H = \PP^1_t \times_{\widehat{g}_1, \widehat{\text{ev}} }
\text{Hom}_{\PP^1_t}( \mc{D}, \PP^1_t \times \Kbm{0,1}{\mb{X},e_2} ).
$$
Again, $(0,[\zeta_C])$ defines a point of $H$ lying over $0\in
\PP^1_t$.  And $\text{pr}_{\PP^1_t}$ is smooth at $(0,[\zeta_C])$ by
base-change. Therefore the image of $\text{pr}_{\PP^1_t}$ is a dense
open subset of $\PP^1_t$.  Let $\ssf{t}\in \PP^1_t$ be a general
point, and define $D\subset \Sigma_1$ to be
$\text{pr}_{\PP^1_t}^{-1}\cap \mc{D}\subset \Sigma_1$.  Then there
exists a morphism $\zeta_D:\mc{D} \rightarrow \Kbm{0,1}{\mb{X},e_2}$
that is a deformation of $\zeta_C$.  If $\ssf{t}$ and $\zeta_D$ are
general, then $\zeta_D$ is a very twisting morphism.  
The curve $D\subset \Sigma_1$ is the image of a section
$\sigma_{1,0}:\PP^1\rightarrow \Sigma_1$.  Denote by
$\widetilde{\zeta}_2:\PP^1 \rightarrow \Kbm{0,1}{\mb{X},e_2}$ the
1-morphism $\widetilde{\zeta} = \zeta_D\circ \sigma_{1,0}$.

\medskip\noindent
By hypothesis, $D\cap
\sigma_{1,1}(\PP^1)$ (resp. $D\cap (\sigma_{1,1}(\PP^1)\cup
\sigma_{1,2}(\PP^1))$) is transverse.
Denote by $\nu:\widetilde{\Sigma}_1 \rightarrow \Sigma_1$ the blowing
up of $\Sigma_1$ at the finitely many intersection points.  Denote by
$\widetilde{\sigma_{1,i}}:\PP^1 \rightarrow \widetilde{\Sigma}_1$ the
strict transform of $\sigma_{1,i}$ for $i=0,1$ (resp. $i=0,1,2$).  In
Case (i), $(D\cdot \sigma_{1,1}(\PP^1))_{\Sigma_1} = (D\cdot D)_{\Sigma_1}
= -\text{deg}(\zeta_1^*\psi_1) > 0$.  If $r=2$, also $(D\cdot
\sigma_{1,2}(\PP^1))_{\Sigma_1} = 0$.  In Case 
(ii), $(D\cdot \sigma_{1,1}(\PP^1))_{\Sigma_1} =
-\text{deg}(\zeta_1^*\psi_1) + \delta > 0$ and $(D \cdot D)_{\Sigma_1}
= -\text{deg}(\zeta_1^*\psi_1) + 2\delta > 0$.  If $r=2$, also
$(D\cdot \sigma_{1,2}(\PP^1))_{\Sigma_1} = \delta$.  In each case, for
the strict transform $\widetilde{D} \subset \widetilde{\Sigma}_1$,
$(\widetilde{D}\cdot \widetilde{D})_{\widetilde{\Sigma}} \geq 0$; it
is precisely $0$ except in Case (ii) if $r=1$.  

\medskip\noindent
Denote by $\widetilde{\pi}_1 = \pi_1\circ \nu$ and denote
$\widetilde{g}_1 = g_1 \circ \nu$.  If $r=1$, denote by
$\widetilde{\zeta}_1:\PP^1 \rightarrow \Kbm{0,r+1}{\mb{X},e_1}$ the
1-morphism associated to the datum
$((\widetilde{\pi}_1:\widetilde{\Sigma}_1 \rightarrow
\PP^1,\widetilde{\sigma}_{1,0},
\widetilde{\sigma}_{1,1}),\widetilde{g}_1)$.  If $r=2$, denote by
$\widetilde{\zeta}_1:\PP^1 \rightarrow \Kbm{0,r+1}{\mb{X},e_1}$ the
1-morphism associated to the datum
$((\widetilde{\pi}_1:\widetilde{\Sigma}_1 \rightarrow \PP^1,
\widetilde{\sigma}_{1,0}, \widetilde{\sigma}_{1,1},
\widetilde{\sigma}_{1,2}), \widetilde{g}_1)$.  In each case, observe
that $\text{ev}_0\circ \widetilde{\zeta}_1 = \text{ev}\circ
\widetilde{\zeta}_2$.  

\medskip\noindent
The following terminology is from ~\cite{BM}.
Denote by $\tau$ the following genus $0$ stable $A$-graph.  There are two
vertices $v_1$ and $v_2$ of degree $e_1$ and $e_2$ respectively.
If $r=1$, there are two flags, $f_0$ and $f_1$, attached to $v_1$.  If
$r=2$, there are three flags, $f_0$, $f_1$ and $f_2$, attached to
$v_1$.  There is one flag attached to $v_2$ and together with $f_0$ it
forms an edge connecting $v_1$ to $v_2$.  If $r=1$, then $f_1$ is a
tail attached to $v_1$.  If $r=2$, both $f_1$ and $f_2$ are tails
attached to $v_1$.  Denote by $\Kbm{}{\mb{X},\tau}$ the associated
Behrend-Manin stack.  In a natural way the pair
$(\widetilde{\zeta}_1,\widetilde{\zeta}_2)$ determines a 1-morphism
$\zeta:\PP^1 \rightarrow \Kbm{}{\mb{X},\tau}$: $\widetilde{\zeta}_1$
corresponds to $v_1$, $\sigma_{1,i}$ corresponds to $f_i$ for each
$i$, and $\widetilde{\zeta}_2$ corresponds to $v_2$.  

\medskip\noindent
There is a canonical contraction of $\tau$ to a stable $A$-graph with
a single vertex of degree $e_1+e_2$ and $r$ tails.  Denote by
$\iota:\Kbm{}{\mb{X},\tau} \rightarrow
\Kbm{0,r}{\mb{X},e_1+e_2}$ the associated 1-morphism.  There is a
canonical combinatorial morphism that is the inclusion of the maximal
subgraph of $\tau$ that contains only the vertex $v_1$.  Denote by
$\alpha:\Kbm{}{\mb{X},\tau} \rightarrow \Kbm{0,r+1}{\mb{X},e_1}$ the
associated 1-morphism.  Finally, denote by
$\beta:\Kbm{0,r+1}{\mb{X},e_1} \rightarrow \Kbm{0,r}{\mb{X},e_1}$ the
isogeny obtained by removing the tail $f_0$.  Observe that
$\text{ev}\circ \beta\circ \alpha$ is 2-equivalent to $\text{ev} \circ
\iota$.  

\medskip\noindent
Of course $\alpha\circ
\zeta = \widetilde{\zeta}_1$ and $\beta\circ \alpha \circ \zeta =
\zeta_1$.  Denote by $\iota \circ \zeta$ by $\zeta$ as well.  The
hypotheses that $\zeta_1$ and $\zeta_2$ are $\text{ev}$-relatively
free imply that $\iota$ is unramified and is \'etale locally a regular
embedding of codimension $1$ at every point of $\zeta(\PP^1)$.  Denote
by $\zeta^*N_\iota$ the locally free sheaf that is the dual of the
pullback of the conormal
sheaf of $\iota$.  
The hypotheses on $\zeta_1$ and $\zeta_2$ imply that
$\alpha$ is unobstructed at every point of $\zeta(\PP^1)$ and $\beta$
is unobstructed at every point of $\alpha(\zeta(\PP^1))$.  Denote by
$\zeta^*T_\alpha$ the locally free sheaf that is the dual of the
pullback of the sheaf of relative differentials of $\alpha$.  Denote
by $\zeta^*\alpha^* T_\beta$ the locally free sheaf that is the dual
of the pullback of the sheaf of relative differentials of $\beta$.
Finally, $\text{ev}$ is unobstructed at every point of
$\beta(\alpha(\zeta(\PP^1)))$ because $\zeta_1$ is
$\text{ev}$-relatively free.  Denote by $\zeta^*\alpha^*\beta^*
T_{\text{ev}}$ the locally free sheaf that is the dual of the pullback
of the sheaf of relative differentials of $\text{ev}$.  
It follows
that also $\text{ev}:\Kbm{0,r}{\mb{X},e_1+e_1} \rightarrow \mb{X}^r$
is unobstructed at every point of $\zeta(\PP^1)$.  Denote by $\zeta^*
T_{\text{ev}}$ the locally free sheaf that is the dual of the pullback
of the sheaf of relative differentials of $\text{ev}$.  
 
\medskip\noindent
There is a filtration of $\zeta^* T_\text{ev}$,
$$
\zeta^* T_\text{ev} = F^0 \supset F^1 \supset F^2 \supset F^3 \supset
F^4=(0),
$$
with associated graded sheaves $F^0/F^1 = \zeta^* N_\iota$, 
$F^1/F^2 = \zeta^* \alpha^*
\beta^* T_\text{ev} = \zeta_1^* T_\text{ev}$, $F^2/F^3 = \zeta^*
\alpha^* T_\beta$, 
and $F^3/F^4 = \zeta^* T_\alpha = \widetilde{\zeta}_2^* T_\text{ev}$.

\medskip\noindent
First, $\zeta^* N_\iota \cong
\OO_{\PP^1}(-\widetilde{\zeta}_1^*\psi_0) \otimes
\OO_{\PP^1}(-\widetilde{\zeta}_2^* \psi)$.  As seen above,
$-\text{deg}(\widetilde{\zeta}_1^* \psi_0) = (\widetilde{D}\cdot
\widetilde{D})_{\widetilde{\Sigma}_1} \geq 0$.  And
$-\text{deg}(\widetilde{\zeta}_2^* \psi) \geq -\text{deg}(\zeta_2^*
\psi) > 0$.  Therefore $\zeta^* N_\iota$ is an ample invertible sheaf.
Next $\zeta_1^* T_\text{ev}$ is ample because $\zeta_1$ is
$\text{ev}$-relatively very free.  Next, $\zeta^* \alpha^* T_\beta
\cong \sigma_{1,0}^* \OO_{\Sigma_1}(\sigma_{1,0}(\PP^1))$ is an ample
invertible sheaf of degree $(D\cdot D)_\Sigma > 0$.  Therefore
$\zeta^* \alpha^* T_\beta$ is an ample invertible sheaf.  
Finally, $\zeta^* T_\alpha = \widetilde{\zeta}_2^* T_\text{ev}$ is
ample because $\widetilde{\zeta}_2$ is very twisting.  Because each of
the associated graded sheaves is ample, also $\zeta^* T_\text{ev}$ is ample.

\medskip\noindent
Also $\text{ev}\circ \zeta = \text{ev} \circ \zeta_1$ is free because
$\zeta_1$ is $\text{ev}$-relatively very free.  Therefore $\zeta$ is
$\text{ev}$-relatively very free.  Of course $\zeta$ maps every point
of $\PP^1$ to a stable map with reducible domain.  Because $\zeta_1$
and $\zeta_2$ are $\text{ev}$-relatively free, each of these stable
maps deforms to a stable map with irreducible domain.  So by
Remark~\ref{rmk-relfree}, a small deformation of $\zeta$ maps a
general point of $\PP^1$ to a stable map with irreducible domain.
Because $\text{ev}_1\circ \zeta=\text{ev}_1 \circ \zeta_1$, it is
nonconstant.  

\medskip\noindent
Finally, if $r=1$ and $\zeta_1$ is very twisting, this is Case (i).
By construction, the self-intersection of
$\widetilde{\sigma}_{1,1}(\PP^1) \subset \widetilde{\Sigma}_1$ is
$0$.  Therefore $\zeta$ is very twisting and
$-\text{deg}(\zeta^*\psi)=0$.  

\end{proof}

\medskip\noindent
A similar argument proves the following result, which shows that if
there exists a twisting morphism $\PP^1 \rightarrow
\Kbm{0,2}{\mb{X},e}$, then for every $r\geq 2$ and every $e'$
sufficiently divisible, there exists a twisting morphism $\PP^1
\rightarrow \Kbm{0,r}{\mb{X},e'}$.  

\begin{prop} \label{prop-operad}
\marpar{prop-operad}
Assume that $\Kbm{0,0}{\mb{X},\epsilon}$ is irreducible for all $\epsilon$
sufficiently divisible.  Let $r \geq 0$, $s_0 \geq 0$, $s_1,\dots,s_r
>0$, and $e_0,e_1,\dots,e_r \geq 0$ be integers.
Let $\zeta_0:\PP^1 \rightarrow \Kbm{0,s_0}{\mb{X},e_0}$ and
$\zeta_i:\PP^1 \rightarrow \Kbm{0,s_i}{\mb{X},e_i}$ be twisting
morphisms for $i=1,\dots,r$.  Then there exists a
twisting morphism $\zeta:\PP^1 \rightarrow \Kbm{0,s}{\mb{X},e}$ where
$s=s_0 + s_1 + \dots + s_{r} - r$ and $e=e_0 + e_1 + \dots + e_r$.   
\end{prop}

\medskip\noindent
Sometimes existence of a very twisting morphism $\zeta:\PP^1
\rightarrow \Kbm{0,1}{\mb{X},e}$ implies existence of a twisting
morphism $\zeta':\PP^1 \rightarrow \Kbm{0,2}{\mb{X},e'}$.  This is
explained in the following proposition.

\begin{notat} \label{notat-rotate}
\marpar{notat-rotate}
Let $((\pi:\Sigma \rightarrow
\PP^1, \sigma:\PP^1 \rightarrow \Sigma),g:\Sigma\rightarrow \mb{X})$
be a datum giving rise to a $1$-morphism $\zeta:\PP^1 \rightarrow
\Kbm{0,1}{\mb{X},e}$.  Assume that the geometric generic fiber of
$\pi$ is irreducible, that $(\sigma(\PP^1)\cdot \sigma(\PP^1))_\Sigma
= 0$, that $g\circ \sigma:\PP^1 \rightarrow \mb{X}$ is very free,
and that $g$ is unramified and \'etale locally a regular
embedding.  Denote by $N$ the locally free sheaf on $\Sigma$ that is
the dual of the kernel of $dg^\dagger:g^* \Omega_\mb{X} \rightarrow
\Omega_\Sigma$.  The linear system $|\sigma(\PP^1)|$ is a
base-point-free pencil.  Denote by $\rho:\Sigma \rightarrow \PP^1$ the
associated morphism.  
Denote by $\tau_1,
\tau_2:\PP^1 \rightarrow \Sigma$ sections of $\rho$ whose images
are general fibers of $\pi$.  Denote by $\zeta':\PP^1 \rightarrow
\Kbm{0,2}{\mb{X},e'}$ the 1-morphism associated to the datum
$((\rho:\Sigma \rightarrow \PP^1,\tau_1,\tau_2),g)$.
\end{notat}

\begin{defn} \label{defn-rotate}
\marpar{defn-rotate}
If $\zeta$ satisfies the hypotheses of Notation~\ref{notat-rotate},
define a \emph{rotation of $\zeta$} to be the 1-morphism $\zeta'$
associated to any choice of $\tau_1$ and $\tau_2$.  
\end{defn}

\begin{prop} \label{prop-rotate}
\marpar{prop-rotate}
If $\zeta$ is very twisting and satisfies the hypotheses of
Notation~\ref{notat-rotate}, then every rotation of $\zeta$ is
twisting. 
\end{prop}

\begin{proof}
The object $L^\vee_{(\pi,g)}$ in the derived category of quasi-coherent
sheaves on $\Sigma$ is represented by the complex,
$$
\begin{CD}
& 0 & & 1 \\
L^\vee_{(\pi,g)}: & \pi^* T_{\PP^1} @>>> \mc{N}
\end{CD}
$$
Denote by $\OO_\Sigma(-D)$ the invertible sheaf
$\OO_\Sigma(-\sigma(\PP^1) -\tau_1(\PP^1) -\tau_2(\PP^1))$.
The object $L^\vee_{(\pi,g)}\otimes \OO_\Sigma(-D)$
fits into a distinguished triangle,
$$
\begin{CD}
\mc{N}(-D)[-1] @>>> L^\vee_{(\pi,g)}(-D) @>>> \pi^* T_{\PP^1}(-D)[0]
@>>> \mc{N}(-D)[0]
\end{CD}
$$
Now $\pi^* T_{\PP^1}(-D)$ has vanishing $h^0$, $h^1$ and $h^2$,
because $R^i\pi_*\OO_\Sigma(-D)$ is $(0)$ for $i=1,2$.  Therefore the
hypercohomology of $L^\vee_{(\pi,g)}(-D)$ equals the hypercohomology
of $\mc{N}(-D)[-1]$.  Because $\zeta$ is very twisting, 
$\mathbf{h}^i(\Sigma,\mc{N}(-D)[-1]) = \mb{h}^i(\Sigma,
L^\vee_{(\pi,g)}(-D)) = 0$ for $i\geq 2$.

\medskip\noindent
The object $L^\vee_(\rho,g)$ in the derived category of quasi-coherent
sheaves on $\Sigma$ is represented by the complex,
$$
\begin{CD}
& 0 & & 1 \\
L^\vee_{(\rho,g)}: & \rho^* T_{\PP^1} @>>> \mc{N} 
\end{CD}
$$
And $\rho^* T_{\PP^1}(-D)$ has vanishing $h^0$, $h^1$ and $h^2$:
$\rho_*$ is zero, and $R^1\rho^*$ is $\OO_{\PP^1}(-1)$.  Therefore the
hypercohomology of $L^\vee_{(\rho,g)}(-D)$ equals the hypercohomology
of $\mc{N}(-D)[1]$.  In particular, $\mb{h}^i(\Sigma,
L^\vee_{(\rho,g)}(-D)) = \mb{h}^i(\Sigma, \mc{N}(-D)[-1]) = 0$ for
$i\geq 2$.

\medskip\noindent
Consider the object $C = R\rho_*
L^\vee_{(\pi,g)}(-\tau_1(\PP^1)-\tau_2(\PP^1))$.  For all $i\geq 3$,
the cohomology sheaf $\mc{H}^i C$ is $(0)$, i.e., $C$ is
quasi-isomorphic to a complex concentrated in degrees $\leq 2$.
Moreover, because this is a family of stable maps, for all $i\leq 0$,
the cohomology sheaf $\mc{H}^i C$ is $(0)$, i.e., $C$ is
quasi-isomorphic to a complex concentrated in degrees $[1,2]$.

\medskip\noindent
Because $g\circ\sigma$ is very free, $\zeta'$ maps the geometric
generic point of $\PP^1$ to an unobstructed point of $\text{ev}$.
Therefore $\mc{H}^2 C$ is a torsion sheaf, which is just
$\mc{H}^2 R\rho_* L^\vee_{(\pi,g)}(-D) \otimes
\OO_{\PP^1}(1)$. Because it is torsion, it is isomorphic to
$\mc{H}^2 C \otimes \OO_{\PP^1}(-1) = \mc{H}^2 R\rho_*
L^\vee_{(\pi,g)}(-D)$.  There is a Leray spectral 
sequence computing the hypercohomology of $L^\vee_{(\pi,g)}(-D)$ in
terms of the hypercohomology of $R\rho_* L^\vee_{(\pi,g)}(-D)$.  
In the spectral sequence, for every differential whose domain is
$H^0(\PP^1, \mc{H}^2 R\rho_* L^\vee_{(\pi,g)}(-D))$, the target vector
space is $(0)$.
Thus there is an injective homomorphism, 
$$
H^0(\PP^1, \mc{H}^2 R\rho_* L^\vee_{(\pi,g)}(-D)) \rightarrow
\mathbf{H}^2(\Sigma, L^\vee_{(\pi,g)}(-D)).
$$
Since $\mathbf{h}^2(\Sigma,L^\vee_{(\pi,g)}(-D))=0$, the torsion sheaf
$\mc{H}^2 R\rho_* L^\vee_{(\pi,g)}(-D) =
(0)$.  So also $\mc{H}^2 C = (0)$, which proves that $\text{ev}$ is
unobstructed at every geomeric point of $\zeta'(\PP^1)$.  This is (i)
of Definition~\ref{defn-twist}.

\medskip\noindent
Moreover, $C$ is quasi-isomorphic to the complex concentrated in
degree $1$, namely $(\zeta')^* T_{\text{ev}}[-1]$.
To prove that $(\zeta')^* T_{\text{ev}}$ is generated by global
sections, it suffices to prove that 
$$
h^1(\PP^1,(\zeta')^* T_{\text{ev}}\otimes \OO_{\PP^1}(-1)) = 0.
$$
This equals $\mathbf{h}^2(\PP^1,C\otimes \OO_{\PP^1}(-1))$.  Because
$C$ is concentrated in a single degree, the Leray spectral sequence
computing hypercohomology of $L^\vee_{(\pi,g)}(-D)$ in terms of
hypercohomology of $C\otimes \OO_{\PP^1}(-1)$ degenerates. In
particular, 
$$
\mathbf{h}^2(\PP^1, C\otimes \OO_{\PP^1}(-1)) =
\mathbf{h}^2(\Sigma,L^\vee_{(\pi,g)}(-D)) = 0.
$$ 
Therefore $(\zeta')^* T_{\text{ev}}$ is generated by global sections.
This is (ii) of Definition~\ref{defn-twist}

\medskip\noindent
Finally, $(\tau_i(\PP^1)\cdot \tau_i(\PP^1))_\Sigma = 0$ for $i=1,2$,
because $\tau_i(\PP^1)$ is a fiber of $\pi$.  Therefore $\zeta^*
\psi_i$ has degree $0$ for $i=1,2$.  This is (iii) of
Definition~\ref{defn-twist}.   
\end{proof}

\section{A very free family of pointed conics} \label{sec-conics}
\marpar{sec-conics}

\begin{hyp} \label{hyp-char0}
\marpar{hyp-char0}
In this section, it is assumed that $\text{char}(K)=0$.
\end{hyp}

\medskip\noindent
Let $\mb{X} \subset \PP^n$ be a general hypersurface of degree $d\leq
n-2$.  The following terminology is from ~\cite{BM}.  
Denote by $\tau$ the stable $A$-graph that has two degree $1$
vertices $v_1$, $v_2$, and edge joining $v_1$ to $v_2$, and two tails
$f_1$, $f_2$ attached to $v_1$ and $v_2$ respectively.  Denote by
$\Kbm{}{\mb{X},\tau}$ the associated Behrend-Manin stack, which is in
fact a scheme.  This stack parametrizes data
$((L_1,x_1,y),(L_2,x_2,y))$ consisting of lines $L_1, L_2 \subset
\mb{X}$ a point $y\in L_1\cap L_2$, and points $x_i\in L_i$ (possibly
equal to $y$).  There is an evaluation morphism
$\text{ev}:\Kbm{}{\mb{X},\tau} \rightarrow \mb{X}^2$ sending a datum
to $(x_1,x_2)$.

\begin{prop}\label{prop-conics1}
\marpar{prop-conics1}
The scheme $\Kbm{}{\mb{X},\tau}$ is an integral, projective scheme of
dimension $3n-2d-1$.
Every irreducible component of the singular locus
$\Kbm{}{\mb{X},\tau}_{\text{sing}}$ has dimension $\leq 2n-d-1$.  
\end{prop}

\begin{proof}
Denote by $\tau_0$ the graph obtained from $\tau$ by removing the
tails $f_1$ and $f_2$.  Then $\Kbm{}{\mb{X},\tau_0} =
\Kbm{0,1}{\mb{X},1} \times_{\mb{X}} \Kbm{0,1}{\mb{X},1}$.  By
~\cite[Theorem V.4.3]{K}, $\Kbm{0,1}{\mb{X},1}$ is smooth, projective and
irreducible of dimension $2n-d-2$.  By ~\cite[Thm. 2.1]{HRS2},
$\text{ev}:\Kbm{0,1}{\mb{X},1} \rightarrow \mb{X}$ is flat of relative
dimension $n-d-1$, and every geometric fiber is connected.  
By straightforward computation, the singular locus of $\text{ev}$
has codimension $\geq n-d$.  Therefore
$\text{pr}_1:\Kbm{}{\mb{X},\tau_0} \rightarrow \Kbm{0,1}{\mb{X},1}$ is
flat of relative dimension $n-d-1$, every geometric fiber is
connected, and the singular locus of $\text{pr}_1$ has codimension
$\geq n-d$ (being the preimage under the flat morphism 
$\text{pr}_2$ of the singular locus of $\text{ev}$).  Thus
$\Kbm{}{\mb{X},\tau_0}$ is an integral, projective scheme of dimension
$3n-2d-3$, and the singular locus has dimension $\leq 2n-d-3$.  The
morphism $\Kbm{}{\mb{X},\tau} \rightarrow \Kbm{}{\mb{X},\tau_0}$ is
the fiber product of two $\PP^1$-bundles. 
\end{proof}

\begin{rmk} \label{rmk-conics1}
\marpar{rmk-conics}
If $d=1$ or $d=2$, then $\Kbm{}{\mb{X},\tau}$ is smooth.
\end{rmk}

\begin{prop} \label{prop-conics2}
\marpar{prop-conics2}
If $n\geq d^2$ and $d\geq 2$, the
geometric generic fiber of $\text{ev}:\Kbm{}{\mb{X},\tau} \rightarrow
\mb{X}\times \mb{X}$ is smooth, connected, nonempty and rationally
connected of dimension $n+1-2d$. 
\end{prop}  

\begin{proof}
Let $(x_1,x_2)\in \mb{X}\times \mb{X}$ be a pair contained in no
common line
contained in $\mb{X}$.  Consider the morphism
$\text{ev}_y:\text{ev}^{-1}(x_1,x_2) \rightarrow \mb{X}$ by
$((L,x_1,y),(L_2,x_2,y)) \mapsto y$.  Because $x_1,x_2$ are contained
in no common line contained in $\mb{X}$, this morphism is a closed
immersion.  Denote by $Y_{x_1,x_2} \subset \PP V$ the image closed subscheme.

\medskip\noindent
Let $V$ be a
vector space of dimension $n+1$ such that $\PP^n = \PP V$.  Let
$L_1,L_2 \subset V$ be $1$-dimensional subspaces corresponding to
$x_1$ and $x_2$.  For $i=1,2$ denote by $\mc{E}_1$ the 
rank $2$ locally $\OO_{\PP V}$-module
$E_i = (L_i^\vee\otimes\OO_{\PP V}) \oplus \OO_{\PP V}(1)$. Denote by
$\phi_i: V^\vee \otimes \OO_{\PP V} \rightarrow E_i$ the unique sheaf
homomorphism such that $\text{pr}_1\circ \phi_i:V^\vee \otimes
\OO_{\PP V} \rightarrow L_i^\vee \otimes \OO_{\PP V}$ and
$\text{pr}_2\circ \phi_i:V^\vee \otimes \OO_{\PP V} \rightarrow
\OO_{\PP V}(1)$ are the canonical surjections.  There is an induced
sheaf homomorphism,
$$
\text{Sym}^d \phi_i: \text{Sym}^d(V^\vee) \otimes \OO_{\PP V}
\rightarrow \text{Sym}^d(E_i) \cong \bigoplus_{k=0}^d
(L_i^\vee)^{\otimes (d-k)} \otimes \OO_{\PP V}(k).
$$
For each $k=0,\dots,d$, denote by $\phi_{i,k}$ composition of
$\text{Sym}^d \phi_i$ with projection onto the $k^\text{th}$ direct
summand.  

\medskip\noindent
Let $F\in \text{Sym}^d(V^\vee)$ be a defining equation of $\mb{X}$.
It is straightforward that $Y_{x_1,x_2}$ is the scheme of simultaneous
zeroes of $\phi_{i,k}(F)$ for $i=1, 1\leq k \leq d$ and $i=2, 1\leq k
\leq d-1$.  Because $n+1-2d \geq 1$, $Y_{x_1,x_2}$ is nonempty and
connected, and every irreducible component has dimension $\geq n+1-2d$.  
In particular, this implies that
$\text{ev}:\Kbm{}{\mb{X},\tau}\rightarrow \mb{X}^2$ is surjective.

\medskip\noindent
By Proposition~\ref{prop-conics1} and Remark~\ref{rmk-conics1}, 
the dimension of the singular locus of $\Kbm{}{\mb{X},\tau}$
is less than $2n-2=\text{dim}(\mb{X}^2)$.  Hence 
the morphism $\text{ev}$ is
generically smooth of relative dimension $n+1-2d$.  So for
$(x_1,x_2)\in \mb{X}^2$ a general pair, $Y_{x_1,x_2}$ is a complete
intersection.  By the adjunction formula, the dualizing sheaf is the
restriction to $Y_{x_1,x_2}$ of the invertible sheaf,
$$
\bigwedge^{n+1}(V^\vee) \otimes \OO_{\PP V}(-n-1)
\otimes \bigotimes_{k=1}^d \lt( (L_1^\vee)^{\otimes (d-k)} \otimes
\OO_{\PP V}(k) \rt)  \otimes \bigotimes_{k=1}^{d-1} \lt(
(L_2^\vee)^{\otimes (d-k)} \otimes \OO_{\PP V}(k) \rt).
$$
In other words, $\omega_Y \cong \OO_{Y}(-n-1+d^2)$.  Because $n\geq
d^2$, $\omega_Y^\vee$ is ample, i.e., $Y_{x_1,x_2}$ is a Fano manifold.
By ~\cite{KMM92c}, ~\cite{Ca}, $Y_{x_1,x_2}$ is
rationally connected. 
\end{proof}

\begin{rmk} \label{rmk-conics2}
\marpar{rmk-conics2}
If $n\geq 3$ and $d=1$, the proposition is still true.  In this case,
for any pair of distinct points $x_1,x_2\in \mb{X}$, the morphism
$\text{ev}_y: \text{ev}^{-1}(x_1,x_2) \rightarrow \mb{X}$ is the
blowing up of $\mb{X}$ at $x_1$ and $x_2$.  
\end{rmk}

\begin{prop} \label{prop-conics}
\marpar{prop-conics}
Let $\mb{X}\subset \PP^n$ be a general hypersurface of degree $d$.  If
either $n\geq d^2$ and $d\geq 2$ or $n\geq 3$ and $d=1$, there exists
a $\text{ev}$-relatively very free 1-morphism $\zeta:\PP^1 \rightarrow
\Kbm{0,2}{\mb{X},2}$ mapping a general point of $\PP^1$ to a stable
map with irreducible domain and with $\text{ev}_1\circ \zeta$
nonconstant.  
\end{prop}

\begin{proof}
The scheme $\mb{X}^2$ is smooth, projective and rationally
connected.  Therefore there exists a family of very free rational
curves dominating $\mb{X}^2$, say $f:D\times \PP^1
\rightarrow \mb{X}^2$.  For each geometric point
$\ssf{t}\in D$, denote by
$\text{pr}_\ssf{t}:\Kbm{}{\mb{X},\tau}_\ssf{t} \rightarrow \PP^1$ the
fiber product,
$$
\Kbm{}{\mb{X},\tau}_\ssf{t} = \PP^1
\times_{f_\ssf{t},\mb{X}^2,\text{ev}} \Kbm{}{\mb{X},\tau}.
$$

\medskip\noindent
By Proposition~\ref{prop-conics1} and Remark~\ref{rmk-conics1}, the
image under $\text{ev}$ of the singular locus of $\Kbm{}{\mb{X},\tau}$
has codimension $\geq 2$.  Therefore a general very free rational
curve does not intersect this locus.  By generic smoothness, for a
general point $\ssf{t} \in D$, $\Kbm{}{\mb{X},\tau}_\ssf{t}$ is smooth
(although $\text{pr}_\ssf{t}$ is only generically smooth).  By
Proposition~\ref{prop-conics2} and Remark~\ref{rmk-conics2}, if
$\ssf{t}$ is general then the geometric generic fiber of
$\text{pr}_\ssf{t}$ is rationally connected, i.e., $\text{pr}_\ssf{t}$
is a rationally connected fibration over a curve.  By ~\cite{GHS},
there exists a section $\sigma_0$
of $\text{pr}_\ssf{t}$, and $\sigma_0(\PP^1)$
is contained in the smooth locus of
$\text{pr}_\ssf{t}$.  Denote by $\widehat{\sigma}_0:\PP^1 \rightarrow
\Kbm{}{\mb{X},\tau}$ the composition of $\sigma_0$ with the projection.

\medskip\noindent
Consider the morphism $\iota: \Kbm{}{\mb{X},\tau}
\rightarrow \Kbm{0,2}{\mb{X},2}$.  Every smooth point of $\text{ev}$
is an unobstructed point.  Every unobstructed point of $\text{ev}$ is
in the open subset where $\iota$ is unramified and \'etale locally a
regular embedding of codimension $1$.  In particular, the image of
$\widehat{\sigma}_0$ is in this open set.  Denote by $N_\iota$ the
locally free sheaf on this open set that is the normal bundle of $\iota$.
  
\medskip\noindent
If $(x_1,x_2)\in \mb{X}^2$ is
general, $\text{ev}^{-1}(x_1,x_2) = Y_{x_1,x_2}$ is contained in the
locus where $\iota$ is unramified and \'etale locally a regular
embedding of codimension $1$.  The restriction of the normal bundle of
$\iota$ to $Y_{x_1,x_2}$ is $\OO_{\PP V}(2)$.  So for a very free
curve in $Y_{x_1,x_2}$, not only is the relative tangent bundle
$T_\text{ev}$ ample, also $N_\iota$ is ample.  Form a comb whose
handle is $\sigma_0$, and whose teeth are very free curves in fibers
of $Y_{x_1,x_2}$.  By ~\cite[II.7.10, II.7.11]{K}, after attaching
sufficiently many 
teeth, the comb deforms to a section $\sigma$ of $\text{pr}_\ssf{t}$
such that both $\sigma^* T_{\text{ev}}$ is ample and $\sigma^*
N_\iota$ is ample.  Denote by $\widehat{\sigma}:\PP^1 \rightarrow
\Kbm{}{\mb{X},\tau}$ the induced morphism.

\medskip\noindent
To distinguish it from $\text{ev}:\Kbm{}{\mb{X},\tau}\rightarrow
\mb{X}^2$, denote by $\text{ev}_a:\Kbm{0,2}{\mb{X},2} \rightarrow
\mb{X}^2$ the evaluation morphism.  On the open set where $\iota$ is
unramified and \'etale locally a regular embedding of codimension $1$,
there is an exact sequence of locally free sheaves,
$$
\begin{CD} 
0 @>>> T_{\text{ev}} @>>> \iota^* T_{\text{ev}_a} @>>> N_\iota @>>> 0
\end{CD}
$$
Therefore, $\widehat{\sigma}^*\iota^* T_{\text{ev}_a}$ is ample.  So
$\iota\circ \sigma:\PP^1 \rightarrow \Kbm{0,2}{\mb{X},2}$ is a
$\text{ev}_a$-relatively very free morphism such that
$\text{ev}_a\circ \iota\widehat{\sigma}$ is nonconstant.  Of course a
general point of $\PP^1$ is mapped to a stable map with reducible
domain.  But this stable map deforms to a stable map with irreducible
domain.  So by Remark~\ref{rmk-relfree}, a small deformation $\zeta$
of $\iota\circ \widehat{\sigma}$ is a $\text{ev}_a$-relatively very
free morphism mapping a general point of $\PP^1$ to a stable map with
irreducible domain and such that $\text{ev}_a\circ \zeta$ is nonconstant.
\end{proof}

\section{Minimal twisting morphisms} \label{sec-outline}
\marpar{sec-outline}

\noindent

\begin{defn} \label{defn-mintwist}
\marpar{defn-mintwist}
Let $\mb{X}$ be a quasi-projective variety with smooth locus
$\mb{X}^0$.  
A twisting morphism $\zeta_0:\PP^1 \rightarrow \Kbm{0,1}{\mb{X}^0,1}$ is
a \emph{minimal twisting morphism} if $\text{deg}(\zeta_0^* h) =1$ and
$\text{deg}(\zeta_0^* x) = 2$, cf. Notation~\ref{notat-xh}.  For a base
scheme $D$, a twisting morphism relative to $D$, $\zeta_0:D\times
\PP^1 \rightarrow \Kbm{0,1}{\mb{X}^0,1}$ is a \emph{minimal twisting
  relative to $D$} if the restriction to every geometric point of $D$
is a minimal twisting morphism.

\medskip\noindent
Denote by $N \subset
\text{Hom}(\PP^1,\Kbm{0,1}{\mb{X}^0,1})$ the locally closed subscheme
parametrizing minimal twisting morphisms.
\end{defn}

\begin{rmk} \label{rmk-mintwist}
\marpar{rmk-mintwist}
\begin{enumerate}
\item[(i)]
In fact $N$ is an open subset.
\item[(ii)]
Let $((\pi:\Sigma\rightarrow \PP^1,\sigma),g)$ be the datum associated
to a twisting morphism $\zeta_0$.  Then $\zeta_0$ is a minimal twisting
morphism iff $g:\Sigma \rightarrow \mb{X}$ is a closed immersion whose
image is a smooth quadric surface and $g\circ \sigma:\PP^1 \rightarrow
\mb{X}$ is a closed immersion whose image is a line.
\item[(iii)]  
Unless $\mb{X}$ is ruled by linear spaces over
a non-uniruled variety, a minimal twisting morphism is
minimal in the following sense: $\text{deg}(\zeta_0^* x) \geq 2$ 
and $\text{deg}(\zeta_0^* h) \geq 1$ for every twisting morphism $\zeta_0$,
cf. ~\cite[Rmk. 5.12]{HS2}.  A variety is called
\emph{quadric type} if a minimal twisting family is minimal in this
sense.  
\item[(iv)]
A smooth hypersurface is quadric type iff the degree is
$\geq 2$.
\end{enumerate}
\end{rmk}


\begin{notat} \label{notat-Tev}
\marpar{notat-Tev}
Let $\mb{X}$ be a quasi-projective variety and denote by $\mb{X}^0$
the smooth locus of $\mb{X}$. 
Denote by
$\Kbm{0,1}{\mb{X},1}_{\text{ev}}$ the maximal open subscheme of
$\Kbm{0,1}{\mb{X}^0,1}$ on which the obstruction group of 
$\text{ev}:\Kbm{0,1}{\mb{X}^0,1}\rightarrow \mb{X}^0$ is zero.  In
other words, a pointed line $[L,x]\in \Kbm{0,1}{\mb{X}^0,1}$ is in
$\Kbm{0,1}{\mb{X},1}_{\text{ev}}$ iff $T_{\mb{X}}\otimes \OO_L$ is
generated by global sections.
Denote
by $T_{\text{ev},\mb{X}}$ the locally free sheaf on
$\Kbm{0,1}{\mb{X},1}$ that is the dual of the sheaf of relative
differentials of $\text{ev}$.
\end{notat}

\medskip\noindent
By definition, if $\zeta_0:\PP^1 \rightarrow \Kbm{0,1}{\mb{X},1}$ is a
twisting family then $\zeta_0(\PP^1) \subset
\Kbm{0,1}{\mb{X},1}_{\text{ev}}$ so that $\zeta_0^*
T_{\text{ev},\mb{X}}$ is defined. 
The twisting family $\zeta_0$ is very
twisting iff $\zeta_0^* T_{\text{ev},\mb{X}}$
is ample.  Typically this is not the case. Let $\mb{X}$ be of quadric
type.  Then for every connected
component of $N$ there exist nonnegative integers $a, b$ 
such that $\zeta_0^*
T_{\text{ev}, \mb{X}} \cong \OO_{\PP^1}^a \oplus \OO_{\PP^1}(1)^b$.
The family $\zeta_0$ is very twisting iff $b>0$ and $a=0$.



\begin{notat} \label{notat-Ms}
\marpar{notat-Ms}
Let $\mb{X}$ be a quasi-projective morphism, and denote by $\mb{X}^0$
the smooth locus of $\mb{X}$. 
Let $\ssf{s}\in \PP^1$ be a point, and let
$[L,x] \in \Kbm{0,1}{\mb{X}^0,1}_{\text{ev}}$ be a
point. 
Denote by $M_{\ssf{s},[L,x]}$ the locally closed subscheme
of $\text{Hom}(\PP^1,\Kbm{0,1}{\mb{X}^0,1})$ that
parametrizes minimal twisting families $\zeta_0:\PP^1 \rightarrow
\Kbm{0,1}{\mb{X}^0,1}$ such that $\zeta_0(\ssf{s}) =
[L,x]$.  
Denote by $M$ the locally closed subscheme of $\PP^1 \times
\Kbm{0,1}{\mb{X}^0,1}_\text{ev} \times N$ 
that
parametrizes triples $(\ssf{s},[L,x],\zeta_0)$ such
that $\zeta_0$ is a minimal twisting family with $\zeta_0(\ssf{s}) =
[L,x]$, i.e., $M$ is the graph of the evaluation
morphism $\text{ev}:\PP^1 \times N \rightarrow
\Kbm{0,1}{\mb{X}^0,1}_{\text{ev}}$.    
\end{notat}

\medskip\noindent
For every $\zeta_0 \in M_{\ssf{s},[L,x]}$,
$\zeta_0^* T_{\text{ev},\mb{X}} \otimes \kappa(\ssf{s})$ equals $T$.
For every connected component $M_{\ssf{s},[L,x],i}$ of
$M_{\ssf{s},[L,x]}$, 
there is a pair of 
nonnegative integers $a,b$ such that $\zeta_0^* T_{\text{ev},\mb{X}}
\cong \OO_{\PP^1}^a \oplus \OO_{\PP^1}(1)^b$ for every $\zeta_0$ in
$M_{\ssf{s},[L,x],i}$.    
The subbundle of
$\zeta_0^* T_{\text{ev},\mb{X}}$ spanned by $\OO_{\PP^1}(1)^b$
restricts to a $b$-dimensional subspace of $T$.  This
defines a morphism from $M_{\ssf{s},[L,x],i}$ to the Grassmannian of $T$.

\begin{notat}\label{notat-q}
\marpar{notat-q}
For each connected component $M_{\ssf{s},[L,x],i} \subset
M_{\ssf{s},[L,x]}$, 
denote by $q_{\ssf{s},[L,x],i}:M_{\ssf{s},[L,x],i} \rightarrow
\text{Grass}(b,T)$ the 
morphism defined above.
The morphism $q_{\ssf{s},[L,x],i}$ is called
\emph{spanning} if there is no proper subspace $T'\subset T$ such that
the image of $q_{\ssf{s},[L,x],i}$ is contained in
$\text{Grass}(b,T')$.  Denote by $G(T)$ the disjoint union over all $b$
of $\text{Grass}(b,T)$.
Denote by $q_{\ssf{s},[L,x]}:M_{\ssf{s},[L,x]} \rightarrow G(T)$
the disjoint union over
all connected components $M_{\ssf{s},[L,x],i}$ of the morphism
$q_{\ssf{s},[L,x],i}$. 
The morphism $q_{\ssf{s},[L,x]}$ is \emph{spanning} if at least one 
$q_{\ssf{s},[L,x],i}$ is spanning. 
\end{notat}

\begin{lem} \label{lem-q}
\marpar{lem-q}
Let $\mb{X}$ be a quasi-projective variety of quadric type, and denote
by $\mb{X}^0$ 
the smooth locus of $\mb{X}$. 
\begin{enumerate}
\item[(i)]
The schemes $M$ and $N$ are smooth, and  
the projection morphism $\text{pr}_{1,2}:M \rightarrow
  \PP^1 \times \Kbm{0,1}{\mb{X}^0,1}_\text{ev}$ is smooth.
\item[(ii)] 
For each connected component $M_{\ssf{s},[L,x],i} \subset
M_{\ssf{s},[L,x]}$, the 
morphism $q_{\ssf{s},[L,x],i}$ is spanning iff the pullback morphism,
$$
q_{\ssf{s},[L,x],i}^*: H^0(\text{Grass}(b,T), \OO_{\text{Grass}}(1)) =
T^\vee 
\rightarrow H^0(M_{\ssf{s},[L,x],i},q_{\ssf{s},[L,x],i}^*
\OO_{\text{Grass}}(1)), 
$$
is injective.
\item[(iii)]
There exists an open subset $U\subset \PP^1 \times
  \Kbm{0,1}{\mb{X}^0,1}_\text{ev}$ such that for every
  $(\ssf{s},[L,x]) \in \PP^1 \times
  \Kbm{0,1}{\mb{X}^0,1}_{\text{ev}}$, $q_{\ssf{s},[L,x]}$ is spanning iff
  $(\ssf{s},[L,x]) \in U$.
\end{enumerate}
\end{lem}

\begin{proof}
\textbf{(i):}  By ~\cite[Thm. II.1.7]{K}, the projection morphism
$\text{pr}_{1,2}:M \rightarrow \PP^1 \times
\Kbm{0,1}{\mb{X}^0,1}_{\text{ev}}$ is smooth if,
$$
h^1(\PP^1,\zeta_0^*
T_{\Kbm{0,1}{\mb{X}^0,1}_{\text{ev}}}\otimes \OO_{\PP^1}(-\ssf{s})) = 0,
$$
for every geometric point  $(\ssf{s},[L_\ssf{s},x_\ssf{s}],\zeta_0)
\in M$.
There is a short exact sequence, 
$$
0 \rightarrow \zeta_0^* T_{\text{ev},\mb{X}}\otimes \OO_{\PP^1}(-\ssf{s})
\rightarrow \zeta_0^* T_{\Kbm{0,1}{\mb{X}^0,1}_{\text{ev}}} \otimes
\OO_{\PP^1}(-\ssf{s}) \rightarrow \zeta_0^* \text{ev}^* T_{\mb{X}^0} \otimes
\OO_{\PP^1}(-\ssf{s}) \rightarrow 0.
$$
Because $\text{ev} \circ \zeta_0: \PP^1 \rightarrow \mb{X}^0$
is a twisting line, $h^1$ of the third term is $0$.  Because
$\zeta_0^* T_{\text{ev},\mb{X}}$ is generated by global sections,
$h^1$ of the first term is $0$.  Therefore, by the long exact sequence
in cohomology, $h^1$ of the middle term is $0$.  

\medskip\noindent
The scheme $\mb{X}^0$ is smooth by definition.  The morphism
$\text{ev}: \Kbm{0,1}{\mb{X}^0,1}_\text{ev} \rightarrow \mb{X}^0$ is
smooth by 
definition.  Hence $\Kbm{0,1}{\mb{X}^0,1}_\text{ev}$ is smooth.  Hence
the product $\PP^1 \times \Kbm{0,1}{\mb{X}^0,1}$ is smooth.  Because
$\text{pr}_{1,2}$ is smooth, the scheme $M$ is smooth.  Because $M$ is
the graph of a morphism, projection $\text{pr}_{1,3}:M \rightarrow
\PP^1 \times N$ is an isomorphism.  Because $\PP^1 \times N$ is
smooth, the scheme $N$ is smooth.

\medskip\noindent
\textbf{(ii):}  
This follows from the definition by taking duals of $T$ and $T'$.  

\medskip\noindent
\textbf{(iii):}
It is not difficult to see that the set $U$ is constructible.
Therefore it suffices to prove it is stable under generization.  Let
$(R,\mathfrak{m})$ be a DVR containing $K$, 
and let $\ssf{s}:\SP(R) \rightarrow
\PP^1$ and $[L,x]:\SP(R) \rightarrow \Kbm{0,1}{\mb{X},1}_{\text{ev}}$
be given.  Denote by $\text{pr}_R:M_R \rightarrow  \SP(R)$ the fiber
product, 
$$
M_R = \SP(R)
\times_{(\ssf{s},[L,x]),\PP^1\times\Kbm{0,1}{\mb{X},1},\text{pr}_{1,2}}
M.
$$
Because $\text{pr}_{1,2}$ is smooth, also $\text{pr}_R$ is smooth.
After replacing $R$ by a finite unramified cover by a DVR, assume that
for every connected component of $M_R$, the geometric generic fiber is
also connected.  Let $M_{R,i}$ be a connected component such that
$q_{\ssf{s},[L,x],i}$ is spanning.

\medskip\noindent
There is a pullback homomorphism of $R$-modules,
$$
q_{R}^*:T^\vee \otimes_K R \rightarrow H^0(M_{R,i},q_R^* \OO_{G(T)}(1)).
$$
Because $q_{\ssf{s},[L,x],i}$ is spanning, the kernel of the closed
fiber of $q_{R}^*$ is injective.
Because $\text{pr}_R$ is smooth, in particular it is flat.  Therefore
the target $R$-module is flat, hence torsion-free.  
So the kernel of $q_R^*$ is a
saturated submodule of $T^\vee$.  Because the kernel of the closed
fiber is zero, also the kernel of the generic fiber is zero.
\end{proof}

\begin{prop} \label{prop-q}
\marpar{prop-q}
Let $\mb{X}$ be a quasi-projective variety of quadric type, and denote
by $\mb{X}^0$ 
the smooth locus of $\mb{X}$.  Denote by $U$ the open set from
Lemma~\ref{lem-q}.  
\begin{enumerate}
\item[(i)]
If $U$ is nonempty, then there exists a very twisting family
$\zeta':\PP^1 \rightarrow \Kbm{0,1}{\mb{X},1}$ such that
$\text{ev}\circ \zeta':\PP^1 \rightarrow \mb{X}^0$ is free.
\item[(ii)] 
If $U$ is nonempty, then for a general deformation $\mb{Y}$ of
$\mb{X}$ there 
exists a very twisting family $\zeta':\PP^1 \rightarrow
\Kbm{0,1}{\mb{Y}^0,1}$ such that $\text{ev}\circ \zeta':\PP^1
\rightarrow \mb{Y}^0$ is free.
\end{enumerate}
\end{prop}

\begin{proof}
\textbf{(i):}
The idea is to construct a map
$\zeta$ from 
a \emph{comb} $C$ to $\Kbm{0,1}{\mb{X},1}$ such that $\zeta$ deforms
to a morphism $\zeta':\PP^1\rightarrow \Kbm{0,1}{\mb{X},1}$ that is
very twisting.  

\medskip\noindent
Define $C_0$ to be a copy of $\PP^1$.  By hypothesis,
there exists a minimal twisting morphism $\zeta_0:C_0 \rightarrow
\Kbm{0,1}{\mb{X},1}$ such that the graph of $\zeta_0$ intersects $U$.
There are nonnegative integers $a,b$ such that $\zeta_0^*
T_{\text{ev},\mb{X}} \cong \OO_{C_0}^a \oplus \OO_{C_0}(1)^b$.  If
$a=0$, then $\zeta_0$ is already very twisting.  Therefore assume that
$a>0$.

\medskip\noindent  
Let
$\ssf{s}_1,\dots \ssf{s}_a \in C_0$ be distinct points and denote by
$[L_1,x_1],\dots, [L_a,x_a]$ their images under $\zeta_0$.  Of course
$\zeta_0^{-1}(U)$ is infinite, so there exist points such that each
pair $(\ssf{s}_i,[L_i,x_i])$ is in $U$.  For each $i$, denote $T_i =
T_{\text{ev},\mb{X}}\otimes \kappa([L_i,x_i])$.  And denote by $T'_i
\subset T_i$ the subspace spanned by the image of $\OO_{C_0}(1)^b$
under restriction to the fiber at $\ssf{s}_i$.  For $i=1,\dots,a$, the
quotient vector space $T_i/T'_i$ is canonically isomorphic to
$H^0(C_0,\zeta_0^* T_{\text{ev},\mb{X}}/\OO_{C_0}(1)^b)$.  Denote this
common vector space by $T/T'$.

\medskip\noindent
Because
each of $(\ssf{s}_i,[L_i,x_i])$ is spanning, for each $i=1,\dots,a$
there exists a curve $C_i$ that is a copy of $\PP^1$ and a minimal
twisting morphism $\zeta_i:C_i \rightarrow \Kbm{0,1}{\mb{X},1}$ that
is in $M_{\ssf{s}_i,[L_i,x_i]}$, and an invertible subsheaf $\mc{L}_i
\cong \OO_{C_i}(1) \subset \zeta_i^* T_{\text{ev},\mb{X}}$ such that
the images of $\mc{L}_i\otimes \kappa(\ssf{s}_i)$ in $T/T'$ span
$T/T'$.

\medskip\noindent
Define $(\iota_i:C_i\hookrightarrow C)_{0\leq i\leq a}$ to be the
initial family of morphisms such that for $i=1,\dots,a$,
$\iota_i(\ssf{s}_i) = \iota_0(\ssf{s}_i)$.  Then $C$ is a connected,
proper, nodal curve whose irreducible components are $C_0,
\dots,C_a$.  Moreover, the dual graph of $C$ is a tree, the vertex of
$C_0$ has valence $a$, and for $i=1,\dots,a$ the vertex of $C_i$ has
valence $1$ and is connected only to the vertex of $C_0$.  The
reducible curve $C$ is a \emph{comb}, the curve $C_0$ is the
\emph{handle}, and the curves $C_1,\dots,C_a$ are the \emph{teeth},
cf. ~\cite[Defn. II.7.7]{K}.

\medskip\noindent
By the universal property of $C$, there is a unique morphism
$\zeta:C\rightarrow \Kbm{0,1}{\mb{X}^0,1}_{\text{ev}}$ such that the
restriction to each $C_i$ is $\zeta_i$.  By ~\cite[Lem. 4.5]{HS2},
$\zeta:C\rightarrow \Kbm{0,1}{\mb{X}^0,1}_{\text{ev}}$ is twisting.

\medskip\noindent
There exists a DVR $R$
containing $K$ and
a proper, flat morphism $\rho:\mc{C} \rightarrow \SP(R)$ whose closed
fiber is $C$, whose generic fiber is $\PP^1$, and such that $\mc{C}$
is regular.  Consider the
relative Hom scheme, 
$$
H_R = \text{Hom}_{\SP(R)}(\mc{C},\SP(R) \times
\Kbm{0,1}{\mb{X}^0,1}_{\text{ev}}).
$$
The morphism $\zeta$ determines a point in the closed fiber of $H_R$.
Because the restriction of $\zeta^* T_{\Kbm{0,1}{\mb{X},1}}$ to every
irreducible component of $C$ is generated by global sections, a
\emph{leaf induction argument} proves that $h^1(C,\zeta^*
T_{\Kbm{0,1}{\mb{X},1}}) = 0$.  By ~\cite[Thm. II.1.7]{K}, the
projection morphism $H_R \rightarrow \SP(R)$ is smooth at $[\zeta]$.
Therefore, after replacing $R$ by a finite unramified cover by a DVR,
there exists a section, i.e., there exists a morphism
$\widetilde{\zeta}:\mc{C} \rightarrow
\Kbm{0,1}{\mb{X}^0,1}_{\text{ev}}$.  Denote by $\zeta':\PP^1
\rightarrow \Kbm{0,1}{\mb{X}^0,1}_{\text{ev}}$ the base-change of
$\widetilde{\zeta}$ over the geometric generic point of $\SP(R)$ (also
base change $K$ to the algebraic closure of the fraction field of $R$).

\medskip\noindent
The claim is that $\zeta'$ is very twisting.  By
~\cite[Lem. 4.6]{HS2}, $\zeta'$ is twisting.  Therefore it only
remains to prove that $(\zeta')^* T_{\text{ev},\mb{X}}$ is ample.
Consider the locally free sheaf $\widetilde{\zeta}^*
T_{\text{ev},\mb{X}}$.  Define $\mc{T} \subset \widetilde{\zeta}^*
T_{\text{ev},\mb{X}}(C_1 + \dots + C_a)$ to be the kernel of the
surjective sheaf homomorphism,
$$
\widetilde{\zeta}^* T_{\text{ev},\mb{X}}(C_1+\dots + C_a) \rightarrow
\bigoplus_{i=1}^a \widetilde{\zeta}^*
T_{\text{ev},\mb{X}}(C_1+\dots+C_a)|_{C_i}/ \mc{L}_i \otimes
\OO_{\mc{C}}(C_1+\dots + C_a)|_{C_i}.
$$
In other words $\mc{T}$ is the \emph{elementary transform up} of
$\widetilde{\zeta}^* T_{\text{ev},\mb{X}}$ along $C_1,\dots, C_a$
determined by $\mc{L}_1, \dots, \mc{L}_a$.  

\medskip\noindent
The sheaf $\mc{T}$ is locally
free, and contains $\widetilde{\zeta}^* T_{\text{ev},\mb{X}}$ as a
subsheaf.  The restriction of $\mc{T}$ to $C_0$ is the sheaf of
meromorphic sections of $\zeta_0^* T_{\text{ev},\mb{X}}$ that have
simple poles in the direction of $\mc{L}_i\otimes \kappa(\ssf{s}_i)$
for each $i=1,\dots,a$.  By construction, this sheaf is isomorphic to
$\OO_{C_0}(1)^{a+b}$.  Also the restriction of $\mc{T}$ to $C_i$ fits
into a short exact sequence,
$$
\begin{CD}
0 @>>> \mc{L}_i @>>> \zeta_i^* T_{\text{ev},\mb{X}} @>>> \mc{T}|_{C_i}
@>>> \mc{L}_i(-\ssf{s}_i) @>>> 0.
\end{CD}
$$
In particular, $\mc{T}|_{C_i}$ is generated by global sections.
Therefore, by
~\cite[Lem 2.11]{HS2}, the sheaf $\mc{T}|_C$ is deformation ample.  By
~\cite[Lem. 2.9]{HS2}, the sheaf $\mc{T}$ on $\mc{C}$ is relatively
deformation ample over $\SP(R)$.  By ~\cite[Lem. 2.8]{HS2}, the
restriction of $\mc{T}$ to the geometric generic fiber of $\mc{C}$ is
ample.  But, of course, the restriction of $\mc{T}$ to the geometric
generic fiber is $(\zeta')^* T_{\text{ev},\mb{X}}$.  Therefore
$\zeta'$ is very twisting.  Also, because $\zeta^* \text{ev}^*
T_{\mb{X}}$ is generated by global sections, also $(\zeta')^*
\text{ev}^* T_{\mb{X}}$ is generated by global sections, i.e.,
$\text{ev}\circ \zeta':\PP^1 \rightarrow \mb{X}^0$ is free.

\medskip\noindent
\textbf{(ii):}
Because $\text{ev}\circ \zeta':\PP^1 \rightarrow \mb{X}^0$ is free, it
is a point at which the relative Kontsevich moduli space is smooth
over the base of the deformation.  Therefore $\text{ev}\circ \zeta'$
deforms to $\mb{Y}$.  By ~\cite[Prop. 4.8]{HS2}, the very twisting
morphism $\zeta'$ also deforms to $\mb{Y}$.
\end{proof}

\medskip\noindent
Because of Proposition~\ref{prop-q}, to prove Theorem~\ref{thm-main}
it suffices to prove that there exists a hypersurface $\mb{X}$, a point
$\ssf{s}$, and a point
$[L,x] \in \Kbm{0,1}{\mb{X}^0,1}$
such that $q_{\ssf{s},[L,x]}$ is spanning.
The boundary case is when $n=d^2$.
This case is the most difficult, and implies the result for all
$n\geq d^2$. 
In this case the subspace of $T$ has
dimension $1$, i.e., $q_{\ssf{s},[L,x]}$ is a morphism
$M_{\ssf{s},[L,x]}E \rightarrow \PP T$.
To prove that $q_{\ssf{s},[L,x]}$ is spanning, it suffices to prove there
exists an element $\zeta_0\in M_{\ssf{s},[L,x]}$ at which
$q_{\ssf{s},[L,x]}$ is 
smooth -- this will even prove that the image of $q_{\ssf{s},[L,x]}$ contains
a dense open subset of $\PP T$.
To see $q_{\ssf{s},[L,x]}$ is smooth at $\zeta_0$, 
it suffices to prove that the 
derivative of $q_{\ssf{s},[L,x]}$ at $\zeta_0$ is surjective.
This is a deformation theory computation that is the
heart of this note.

\medskip\noindent
Before proceeding to this computation, the following proposition shows
that the inequality $n\geq d^2$ in Theorem~\ref{thm-main} is necessary.

\begin{prop} \label{prop-ineq}
\marpar{prop-ineq}
Let $K$ be an algebraically closed field of characteristic $0$.  
Let $(d,n)$ be a pair of
positive integers such that $d \leq \frac{n}{2}$.
Let $\mb{X} \subset \PP^n$ be a hypersurface of degree $d$, and denote
by $\mb{X}^0$ the smooth locus of $\mb{X}$.
If there exists a very twisting
morphism $\zeta:\PP^1 \rightarrow \Kbm{0,1}{\mb{X}^0,1}$, then $n\geq d^2$.
\end{prop}

\begin{proof}
By ~\cite[Defn. 4.3 (ii)]{HS2}, the stable map $\text{ev}\circ \zeta:
\PP^1 \rightarrow \mb{X}^0$ is unobstructed.  This implies that the
relative Kontsevich moduli space of the universal family of degree $d$
hypersurfaces over $\PP \text{Sym}^{d}(K^{n+1})$ is smooth over $\PP
\text{Sym}^d(K^{n+1})$ at $([\mb{X}],[\text{ev}\circ\zeta])$.  
So for a general
deformation of $\mb{X}$, the stable map deforms as well.  By
~\cite[Prop. 4.8]{HS2}, also the very twisting family deforms.
So for a general hypersurface, there is also a very twisting
morphism.  Therefore assume that $\mb{X}$ is general (this is made
precise in the next paragraph).

\medskip\noindent
By ~\cite{HRS2}, if $\mb{X}$ is general then every Kontsevich
moduli space is irreducible.  Also, for every degree $e$ there is a free
rational curve of degree $e$: because $d\leq n$ there is a free line $L$
on $\mb{X}$, 
and a finite degree $e$ morphism $\PP^1 \rightarrow L$ is free.
The locus of free rational curves is open.
Therefore, a general point of the Kontsevich
moduli space parametrizes a free rational curve.  By the same argument
in the last paragraph, $\zeta$ can be deformed so that
$\text{ev}\circ \zeta:\PP^1 \rightarrow \mb{X}$ is a free rational curve.

\medskip\noindent
Consider the Picard group of $\Kbm{0,1}{\PP^n,1}$.  There is the
projection morphism $\text{pr}:\Kbm{0,1}{\PP^n,1} \rightarrow
\Kbm{0,0}{\PP^n,1} = \text{Grass}(2,n+1)$.  And there is the
projection morphism $\text{ev}:\Kbm{0,1}{\PP^n,1} \rightarrow \PP^n$.
Denote by $x$ the
first Chern class of $\text{pr}^* \OO_{\text{Grass}}(1)$, and
denote by $h$ the first Chern class of $\text{ev}^* \OO_{\PP^n}(1)$.
For $\Kbm{0,1}{\PP^n,1}$, $\text{Chow}^1 = \ZZ\{x,h\}$.  There is a
tautological class $\psi = x-2h$.  For any family of pointed lines
parametrized by a base $B$, say $\zeta_{\PP^n} =(\pi:\Sigma
\rightarrow B, \sigma:B\rightarrow \Sigma,g:\Sigma\rightarrow \PP^n)$,
$C_1(\sigma^* \OO_{\Sigma}(\sigma(B))) = -\zeta^*\psi$.  

\medskip\noindent
It is straightforward to compute that,
$$
C_1(T_{\text{ev},\PP^n}) = nx-(n-1)h = (n+1)x+n\psi.
$$  
Denote the universal family of pointed lines by,
$$
\widetilde{\zeta} = (\widetilde{\pi}:\widetilde{\Sigma}
\rightarrow \Kbm{0,1}{\PP^n,1}, \widetilde{\sigma}:\Kbm{0,1}{\PP^n,1}
\rightarrow \widetilde{\Sigma}, \widetilde{g}:\widetilde{\Sigma}
\rightarrow \PP^n).
$$
It is straightforward to compute that,
$$
C_1(\widetilde{\pi}_* (\widetilde{g}^* \OO_{\PP^n}(d) \otimes
\OO_{\widetilde{\Sigma}}(-\text{Image}(\widetilde{\sigma})))) =
\frac{d(d+1)}{2} x -dh = d^2h + \frac{d(d+1)}{2}\psi.
$$
Therefore, if $\zeta_\mb{X}:B \rightarrow
\Kbm{0,1}{\mb{X},1}_{\text{ev}}$ is a morphism and $\zeta_{\PP^n}$ is
the associated morphism to $\Kbm{0,1}{\PP^n,1}$, then,
$$
\begin{array}{c}
C_1(\zeta_{\mb{X}}^* T_{\text{ev},\mb{X}}) = C_1(\zeta_{\PP^n}^*
T_{\text{ev},\PP^n}) - C_1(\pi_* (g^* \OO_{\PP^n}(-d) \otimes
\OO_\Sigma(-\sigma(B)))) = \\ 
(n-\frac{d(d+1)}{2})\zeta^*x - (n-d-1)\zeta^*h = 
(n+1-d^2)\zeta^*h + (n-\frac{d(d+1)}{2})\zeta^*\psi.
\end{array}
$$
In particular, this holds for the very twisting family $\zeta:\PP^1
\rightarrow \Kbm{0,1}{\mb{X},1}$.  

\medskip\noindent
There exists a short exact sequence,
$$
\begin{CD}
0 @>>> \zeta^* T_{\text{ev},\mb{X}} @>>> \zeta^*
T_{\Kbm{0,1}{\mb{X},1}} @>>> \zeta^* \text{ev}^* T_X @>>> 0
\end{CD}
$$
Because $\zeta$ is very twisting, the first term is positive.
Because $\text{ev}\circ \zeta$ is free, the third term is semipositive.
Therefore $\zeta^* T_{\Kbm{0,1}{\mb{X},1}}$ is semipositive.  Because
$\zeta^* \text{pr}^* T_{\Kbm{0,0}{\mb{X},1}}$ is a quotient of
$\zeta^* T_{\Kbm{0,1}{\mb{X},1}}$, it is also semipositive. 
The derivative $d(\text{pr}\circ \zeta):\OO_{\PP^1}(2) \rightarrow
\zeta^* \text{pr}^* T_{\Kbm{0,0}{\mb{X},1}}$ is zero iff
$\text{pr}\circ \zeta$ is constant.  Therefore either $C_1(\zeta^*
T_{\Kbm{0,1}{\mb{X},1}})$ has positive degree or $\text{pr}\circ
\zeta$ is constant.

\medskip\noindent
The scheme $\Kbm{0,1}{\mb{X},1}$ is simply the Fano scheme of lines on
$\mb{X}$.  The canonical bundle of the Fano scheme is straightforward
to compute, giving,
$$
\text{pr}^* C_1(T_{\Kbm{0,0}{\mb{X},1}}) = (n+1 - \frac{d(d+1)}{2}) x.
$$ 
So either $n+1 >
\frac{d(d+1)}{2}$, or else $\text{pr}\circ \zeta$ is constant.
In the second case $\text{deg}(\zeta^* x)=0$, so that,
$$
\text{deg}(C_1(\zeta^* T_{\text{ev},\mb{X}})) =
-(n-d-1)\text{deg}(\zeta^* h).
$$  
Since $d \leq n-1$ this is nonpositive, contradicting that
$\zeta^* T_{\text{ev},\mb{X}}$ is ample.  Therefore $n+1 >
\frac{d(d+1)}{2}$, i.e. $n-\frac{d(d+1)}{2} \geq 0$.

\medskip\noindent
Because $\zeta$ is very twisting,
$\sigma^*\OO_{\Sigma}(\sigma(\PP^1))$ has nonnegative
degree, i.e., $\text{deg}(\zeta^* \psi) \leq 0$.  So
$(n-\frac{d(d+1)}{2})\text{deg}(\zeta^* \psi) \leq 0$.  Therefore,
$$
\text{deg}(C_1(\zeta^* T_{\text{ev},\mb{X}})) \leq
(n+1-d^2)\text{deg}(\zeta^* h).
$$
Because $\zeta^* T_{\text{ev},\mb{X}}$ is ample,
$(n+1-d^2)\text{deg}(\zeta^* h) > 0$.  And $\text{deg}(\zeta^* h) \geq
0$.  Therefore $n+1-d^2 >0$, i.e., $n\geq d^2$.
\end{proof}

\section{Notation} \label{sec-not}
\marpar{sec-not}
\noindent
In this section, the notation for the computation is introduced.  A
specific homogeneous polynomial $G$ of degree $d$ is given, and the
hypersurface 
$\mb{X}$ is $\mb{V}(G)\subset \PP^n$.  A minimal twisting family is specified,
$$
\zeta_0 = (\pi:\Sigma \rightarrow B,\sigma:B\rightarrow \Sigma,
f_0:\Sigma \rightarrow \mb{X}).
$$
In a later section it is proved that $\zeta_0$ is a smooth point of
$\text{Hom}(\PP^1,\Kbm{0,1}{\mb{X},1})$.  An Artin local scheme $D$ is
specified that is a closed subscheme of the first-order neighborhood
of $[\zeta_0]$ in $\text{Hom}(\PP^1,\Kbm{0,1}{\mb{X},1})$.  The family
$\zeta_0$ is ``thickened'' to a family,
$$
\zeta = ((1,\pi):D\times \Sigma \rightarrow D\times B,
(1,\sigma):D\times B \rightarrow D\times \Sigma, f:D\times \Sigma
\rightarrow \mb{X}).
$$
For each closed point $\ssf{s}\in B$, the pointed line
$\zeta_0(\ssf{s})$ is denoted by
$[L_\ssf{s},x_\ssf{s}]$, and the scheme
$M_{\ssf{s},[L_\ssf{s},x_\ssf{s}]}$ 
from Notation~\ref{notat-Ms} is
denoted by $M_\ssf{s}$.  Similarly, the morphism
$q_{\ssf{s},[L_\ssf{s},x_\ssf{s}]}$ from Notation~\ref{notat-q} is
denote by $q_\ssf{s}$. 
The subscheme $D\cap M_{\ssf{s}} \subset D$ is denoted by 
$D_\ssf{s}$.
It turns
out that to prove $dq_{\ssf{s}}|_{\zeta_0}$ is surjective, it suffices to
restrict to the Zariski tangent space $T_0 D_\ssf{s} \subset
T_{[\zeta_0]} M_\ssf{s}$.  This restriction is denoted
$d'q_{\ssf{s}}$.  

\medskip\noindent
Let $K$ be an algebraically closed field with
$\text{char}(K) \neq 2$.  
Let $d \geq 3$ be an integer and denote $n=d^2$.  Denote by $I_d$ the
set with $d^2-4$ elements,
$$
\{(i,j) \in \ZZ^2| 0\leq i,j\leq d-1, (i,j)\neq
(0,0),(0,1),(1,0),(1,1)\}.
$$  
Denote by $V$ the $(n+1)$-dimensional
$K$-vector space with ordered basis,
$$
(\mb{a}_0,\mb{a}_1,\mb{a}_2,\mb{a}_3)\cup(\mb{b}_{(i,j)}|(i,j)\in I_d)
\cup (\mb{c}).
$$
Denote the dual ordered basis of $V^\vee$ by,
$$
(X_0,X_1,X_2,X_3)\cup(Y_{(i,j)}| (i,j) \in I_d) \cup (Z).
$$
Denote by $V_a$ the subspace of $V$ generated by $\mb{a}_0,\dots,\mb{a}_3$.

\medskip\noindent
Given a pair of nonnegative integers, $i$ and $j$, denote by
$k=k(i,j)$ the minimum.  Denote by $G \in \text{Sym}^d(V^)$ the
homogeneous polynomial,
$$
G = (X_0X_3-X_1X_2) X_3^{d-2} + \sum_{(i,j)\in I_d} X_0^k
X_1^{i-k}X_2^{j-k} X_3^{d-1-i-j+k} Y_{(i,j)} + \sum_{(i,j) \in I_d}
Y_{(i,j)}Z \gamma_{(i,j)} + Z^2 \gamma_{z},
$$
where the elements $\gamma_{(i,j)}$ and $\gamma_z$ are homogeneous
polynomials of degree $d-2$ in the variables $X_0,X_1,X_2,X_3$.
Denote by $\mb{X}\subset \PP V$ the hypersurface defined by $G$.

\medskip\noindent
Denote by $\Sigma$ the surface,
$$
\PP^1_s \times \PP^1_t = \{([S_0:S_1],[T_0:T_1]) | [S_0:S_1],[T_0:T_1]
\in \PP^1 \}. 
$$
For each pair of integers $(i,j)$, denote by $\OO_{\Sigma}(i,j)$ the
invertible sheaf $\text{pr}_{\PP^1_s}^*\OO_{\PP^1_s}(i)\otimes
\text{pr}_{\PP^1_t}^* \OO_{\PP^1_t}(j)$.  

\medskip\noindent
Denote by $B$ the curve $\PP^1_s$.  Denote by $\pi:\Sigma \rightarrow
B$ projection onto the first factor. Denote by $\sigma:B\rightarrow
\Sigma$ the section of $\pi$,
$$
\sigma([S_0:S_1]) = ([S_0:S_1],[1:0]).
$$

\medskip\noindent
Denote by $E$ the $K$-vector space with ordered basis,
$$
(u_{(i,j)}^0,u_{(i,j)}^1,u_{(i,j)}^2,u_{(i,j)}^3|(i,j)\in I_d) \cup
(v^0,v^1,v^2,v^3).
$$
Denote by $D'$ the local, Artin $K$-scheme,
$$
D' = \SP \text{Sym}^\bullet(E)/\text{Sym}^2(E)\cdot
\text{Sym}^\bullet(E).
$$

\medskip\noindent
There is a sheaf homomorphism, $V^\vee \otimes_K \OO_{D'} \rightarrow
V_a^\vee \otimes_K \OO_{D'}$ by,
$$
\lt\{
\begin{array}{crcc} 
X_l & \mapsto & X_l, & l=0,1,2,3 \\
\\
Y_{(i,j)} & \mapsto & \sum_{l=0}^3 u_{(i,j)}^l X_l, & (i,j)\in I_d \\
\\
Z & \mapsto & \sum_l v^l X_l
\end{array} \rt.
$$
This induces a morphism of schemes, $\phi':D'\times V_a
\rightarrow D' \times V$.
There is a unique morphism of schemes $\psi':D' \times \Sigma
\rightarrow D' \times \PP V_a$ such that $(\psi')^*\OO_{\PP V_a}(1) =
\OO_{\Sigma}(1,1)$ and such that,
$$
X_0 \mapsto S_0T_0, \ X_1 \mapsto S_0T_1, \ X_2 \mapsto S_1T_0, \ X_3
\mapsto S_1 T_1.
$$
Define $f':D' \times \Sigma \rightarrow D'\times \PP V$ to be the
composition $f' = \phi'\circ \psi'$.

\medskip\noindent
If $(i,j)$ is a pair that is not in $I_d$, define $u_{(i,j)}^l$ to be
$0$.  Denote by $D\subset D'$ the closed subscheme whose ideal
sheaf is,
$$
\mb{I}(D) = \langle u_{(i-1,j-1)}^0 + u_{(i-1,j)}^1 + u_{(i,j-1)}^2 +
u_{(i,j)}^3 | 0\leq i,j \leq d \rangle.
$$
Define $\phi$ to be the restriction of $\phi'$ to $D\times V$, define
$\psi$ to be the restriction of $\psi'$ to $D\times \Sigma$, and
define $f$ to be the composition $f=\phi\circ \psi$.  Define $\phi_0$,
$\psi_0$ and $f_0$ to be the restrictions to the closed point $0\in D$.

\begin{lem} \label{lem-Gis0}
\marpar{lem-Gis0}
The preimage $f^* G$ in
$H^0(D\times \Sigma,f^* \OO_{\PP V}(d))$ equals $0$.  
Therefore the
morphism $f$ factors through $D\times \mb{X} \subset D\times \PP V$.
\end{lem}

\begin{proof}
For every $(i,j)\in I_d$, $f^*Y_{(i,j)} \in \mathfrak{m}\cdot
H^0(D\times \Sigma, f^*\OO_{\PP V}(d))$.  Similarly for $f^*Z$.
Because $\mathfrak{m}^2 =0$, this implies that every term that
involves $Y_{(i,j)}Z$ or $Z^2$ pulls back to $0$.  So these terms are
ignored.  Without these terms
the preimage of $f^* G$ is,
$$
\begin{array}{c}
(S_0T_0\cdot S_1T_1 - S_0T_1\cdot S_1T_0)(S_1T_1)^{d-2} + \\
\sum_{(i,j)\in I_d}S_0^iS_1^{d-1-i}T_0^j T_1^{d-1-j}(u_{(i,j)}^0 S_0T_0
+ u_{(i,j)}^1 S_0T_1 + u_{(i,j)}^2 S_1T_0 + u_{(i,j)}^3 S_1T_1).
\end{array}
$$
The first term is $0$ because $S_0T_0\cdot S_1T_1 = S_0T_1\cdot
S_1T_0$.  Gathering like monomials $S_0^iS_1^{d-i}T_0^jT_1^{d-j}$,
and using the notation that $u_{(i,j)}^l=0$ if $(i,j)\not\in I_d$,
the remaining terms sum to,
$$
\sum_{0\leq i,j \leq d} (u_{(i-1,j-1)}^0 + u_{(i-1,j)}^1 +
u_{(i,j-1)}^2 + u_{(i,j)}^3) S_0^i S_1^{d-i} T_0^j T_1^{d-j}.
$$
By definition, each element $u_{(i-1,j-1)}^0 +u_{(i-1,j)}^1 +
u_{(i,j-1)}^2 + u_{(i,j)}^3$ equals $0$ in $\OO_D$.  Therefore $f^* G$
equals $0$.
\end{proof}

\medskip\noindent
Denote,
$$
\begin{array}{c}
\zeta_{\PP V} = ((1,\pi):D\times \Sigma \rightarrow D\times B,
(1,\sigma):D \times B \rightarrow D\times \Sigma, f:\Sigma \rightarrow
\PP V), \\
\zeta_{\mb{X}} = ((1,\pi):D\times \Sigma \rightarrow D\times B,
(1,\sigma):D \times B \rightarrow D\times \Sigma, f:\Sigma \rightarrow
\mb{X}), \\
\zeta_{\PP V,0} = (\pi:\Sigma \rightarrow B, \sigma:B \rightarrow
\Sigma, f_0: \Sigma \rightarrow \PP V), \\
\zeta_{\mb{X},0} = (\pi:\Sigma \rightarrow B, \sigma:B \rightarrow
\Sigma, f_0: \Sigma \rightarrow \mb{X})
\end{array}
$$

\medskip\noindent
For each $\ssf{s}=[\ssf{s}_0:\ssf{s}_1] \in \PP^1_s$, denote by
$D_\ssf{s} \subset D$ the closed subscheme whose ideal is,
$$
\mb{I}(D_\ssf{s}) = \langle \ssf{s}_0 u_{(i,j)}^0 + \ssf{s}_1
u_{(i,j)}^2, \ssf{s}_0 u_{(i,j)}^1 + \ssf{s}_1 u_{(i,j)}^3 | (i,j) \in
I_d \rangle + \langle \ssf{s}_0 v^0 + \ssf{s}_1 v^2, \ssf{s}_0 v^1 +
\ssf{s}_1 v^3 \rangle.
$$
Denote by $\Lambda_\ssf{s} \subset V_a \subset V$ the subspace,
$$
\text{span}\{ \ssf{s}_0 \mb{a}_0 + \ssf{s}_1 \mb{a}_2, \ssf{s}_0
\mb{a}_1 + \ssf{s}_1 \mb{a}_3 \}.
$$
Denote $L_\ssf{s} = \PP \Lambda_\ssf{s} \subset \PP V$, and denote
$x_\ssf{s} \in L_\ssf{s}$ the point $[\ssf{s}_0 \mb{a}_0 + \ssf{s}_1
\mb{a}_2] \in L_\ssf{s}$.

\begin{lem} \label{lem-Ds}
\marpar{lem-Ds}
The closed subscheme $D_\ssf{s} \subset D$ is the maximal subscheme over which
$f:D_\ssf{s} \times \pi^{-1}(\ssf{s}) \rightarrow D_\ssf{s} \times \PP
V$ factors through $D_\ssf{s} \times L_\ssf{s} \subset D_\ssf{s}
\times \PP V$.
\end{lem}

\begin{proof}
The image of $D\times \pi^{-1}(\ssf{s})$ under $\phi$ is $V_\ssf{s}$,
where $V_\ssf{s} \subset D\times V$ is the vector subbundle with generators,
$$
\begin{array}{c}
\ssf{s}_0 \mb{a}_0 + \ssf{s}_1 \mb{a}_2 + \sum_{(i,j)\in
  I_d}(\ssf{s}_0 u_{(i,j)}^0 + \ssf{s}_1 u_{(i,j)}^2)\mb{b}_{(i,j)} +
  (\ssf{s}_0 v^0 + \ssf{s}_1 v^2) \mb{c}, \\
\ssf{s}_0 \mb{a}_1 + \ssf{s}_1 \mb{a}_3 + \sum_{(i,j)\in
  I_d}(\ssf{s}_0 u_{(i,j)}^1 + \ssf{s}_1 u_{(i,j)}^3)\mb{b}_{(i,j)} +
  (\ssf{s}_0 v^1 + \ssf{s}_1 v^3) \mb{c}
\end{array}
$$
Therefore the maximal closed subscheme $D_\ssf{s}$ such that
$D_\ssf{s} \times_D V_\ssf{s}$ is contained in $D_\ssf{s} \times
\Lambda_\ssf{s}$ is the closed subscheme where the coefficients of
$\mb{b}_{(i,j)}$ and $\mb{c}$ in the generators equal $0$, i.e., the
closed subscheme $D_\ssf{s}$ defined above.
\end{proof}

\section{Description of some coherent sheaves} \label{sec-coh}
\marpar{sec-coh}

\noindent
Denote by $\text{ev}:\Kbm{0,1}{\PP V, 1} \rightarrow \PP V$ the
evaluation morphism.  Denote by $T_{\text{ev},\PP V}$ the locally free
sheaf on $\Kbm{0,1}{\PP V,1}$ that is the dual of the sheaf of
relative differentials of $\text{ev}$.  Denote by
$\text{ev}:\Kbm{0,1}{\mb{X},1} \rightarrow \mb{X}$ the evaluation
morphism and denote by $\Kbm{0,1}{\mb{X},1}_{\text{ev}} \subset
\Kbm{0,1}{\mb{X},1}$ the open subscheme where $\text{ev}$ is smooth.
Denote by $T_{\text{ev},\mb{X}}$ the locally free sheaf on
$\Kbm{0,1}{\mb{X},1}_{\text{ev}}$ that is the dual of the sheaf of
relative differentials of $\text{ev}$.  

\medskip\noindent
The families $\zeta_{\PP V}$ and $\zeta_{\mb{X}}$ in the last section
define morphisms $\zeta_{\PP V}: D\times B \rightarrow \Kbm{0,1}{\PP
  V,1}$ and $\zeta_{\mb{X}}:D \times B \rightarrow \Kbm{0,1}{\mb{X},
  1}$.  In this section it is proved that the image of
$\zeta_{\mb{X}}$ is contained in $\Kbm{0,1}{\mb{X},1}_{\text{ev}}$,
and explicit descriptions are given for $\zeta_{\PP V}^*
T_{\text{ev},\PP V}$ and $\zeta_{\mb{X}}^* T_{\text{ev},\mb{X}}$.
This description will be the basis for the computation of $dq_\ssf{s}$
in the next section.

\medskip\noindent
Denote by $\widetilde{\psi}:D\times \Sigma \rightarrow D \times B
\times \PP V_a$ the unique morphism such that
$\text{pr}_{1,3}\circ \widetilde{\psi} = \psi$ and $\text{pr}_{2}\circ
\widetilde{\psi} = \pi \circ \text{pr}_\Sigma$.  Similarly, denote by
$\widetilde{f}: D\times \Sigma \rightarrow D\times B \times \PP V$ the
unique morphism such that $\text{pr}_{1,3}\circ \widetilde{f} = f$ and
$\text{pr}_2 \circ \widetilde{f} = \pi \circ \text{pr}_\Sigma$.
Denote by $\widetilde{\psi}:D\times B \times \PP V_a \rightarrow
D\times B \times \PP V$ the unique morphism such that $\text{pr}_{1,3}
\circ \widetilde{\psi} = \psi \circ \text{pr}_{1,3}$ and $\text{pr}_2
\circ \widetilde{\psi} = \text{pr}_2$.

\medskip\noindent
Denote by $N_{\widetilde{\phi}}$ the normal bundle of the closed
immersion $\widetilde{\phi}$. 
Denote by $N_{\widetilde{f}}$ 
the normal bundle of the closed immersion $\widetilde{f}$.  There
is a short exact sequence of locally free sheaves on $D\times \Sigma$,
$$
0 \rightarrow N_{\widetilde{\phi}} \xrightarrow{d \widetilde{\psi}} 
N_{\widetilde{f}} \rightarrow
V\otimes_K f^*\OO_{\PP V}(1)/\psi(V_a)\cdot f^*\OO_{\PP
  V}(1) \rightarrow 0.
$$
Using the Euler exact sequence for the tangent bundle of $\PP V$ and
the tangent bundle of $\PP^1_t$, 
the locally free sheaf $N_{\widetilde{\phi}}$ is the cokernel of the
sheaf homomorphism,
$$
\text{pr}_{\Sigma}^* \OO_{\Sigma}(0,1) \{ t_0^\vee, t_1^\vee \} 
\rightarrow \text{pr}_{\Sigma}^*\OO_{\Sigma}(1,1)\{
\mb{a}_0, \mb{a}_1, \mb{a}_2, \mb{a}_3 \},
$$
defined by,
$$
t_0^\vee \mapsto s_0 \mb{a}_0 + s_1 \mb{a}_2, \ t_1^\vee \mapsto s_0
\mb{a}_1 + s_1 \mb{a}_3.
$$
The cokernel is the sheaf homomorphism,
$$
\text{pr}_{\Sigma}^* \OO_{\Sigma}(1,1)\{ \mb{a}_0, \mb{a}_1, \mb{a}_2,
\mb{a}_3 \} \rightarrow
\text{pr}_{\Sigma}^* \OO_{\Sigma}(2,1) \{ \mb{e}_0, \mb{e}_1 \},
$$
defined by,
$$
\mb{a}_0 \mapsto s_1 \mb{e}_0, \ \mb{a}_1 \mapsto
s_1 \mb{e}_1, \ \mb{a}_2 \mapsto -s_0 \mb{e}_0, \ \mb{a}_3 \mapsto
-s_0 \mb{e}_1. 
$$
This cokernel is identifed with $N_{\widetilde{\phi}}$.

\medskip\noindent
There is a canonical surjective sheaf homomorphism $f^* T_{\PP V}
\rightarrow N_{\widetilde{f}}$.  Moreover, by the Euler exact sequence
for $\PP V$, there is a canonical sheaf homomorphism $V\otimes_K
\OO_{\PP V}(1) \rightarrow T_{\PP V}$.  Together, these give a
canonical sheaf homomorphism $V\otimes_K \text{pr}_\Sigma
\OO_{\Sigma}(1,1) \rightarrow N_{\widetilde{f}}$.  For each element
$\mb{v}$ of $V$, also denote by $\mb{v}$ the induced sheaf
homomorphism $\text{pr}_\Sigma^*
\OO_{\Sigma}(1,1) \rightarrow N_{\widetilde{f}}$, and
the induced global section of $N_{\widetilde{f}}\otimes \text{pr}_\Sigma^*
\OO_{\Sigma}(-1,-1)$.  

\medskip\noindent
In particular there are global sections of $N_{\widetilde{f}} \otimes
\text{pr}_\Sigma^* \OO_{\Sigma}(-1,-1)$:
$\mb{c}$ and $\mb{b}_{(i,j)}$, for $(i,j)\in I_d$. 
The images of these global sections in $V\otimes_K
\OO_{D\times \Sigma}/\widetilde{f}(V_a)$ give a basis as a free
$\OO_{D\times \Sigma}$-module.  With the
isomorphism from the last paragraph this defines an isomorphism,
$$
N_{\widetilde{f}} \cong \text{pr}_{\Sigma}^* \OO_{\Sigma}(2,1)\{ \mb{e}_0,
\mb{e}_2 \} \oplus \text{pr}_{\Sigma}^* \OO_{\Sigma}(1,1)\{
\mb{b}_{(i,j)} | (i,j) \in I_d \} \oplus \text{pr}_{\Sigma}^* 
\OO_{\Sigma}(1,1) \{ \mb{c} \}.
$$



\medskip\noindent
The locally free sheaf $\zeta_{\PP V}^* T_{\text{ev},\PP V}$ is
canonically isomorphic to $(1\times \pi)_*
N_{\widetilde{f}}(-D\times\sigma(B))$, cf. ~\cite[Sec. 3]{HS2}.
Therefore, 
$$
\zeta_{\PP V}^* T_{\text{ev},\PP V} \cong 
\text{pr}_{B}^* \OO_{B}(2)\{ \mb{e}_0,
\mb{e}_2 \} \oplus \text{pr}_{B}^* \OO_{B}(1)\{ \mb{b}_{(i,j)} | (i,j)
\in I_d \} \oplus \text{pr}_{B}^* 
\OO_{B}(1) \{ \mb{c} \}.
$$
The morphism $T_{\PP_V}|_{\mb{X}} \rightarrow N_{\mb{X}/\PP V}$
induces a sheaf homomorphism,
$$
N_f(-D\times \sigma(B)) \rightarrow f^*(\OO_{\PP V}(d))(-D\times
\sigma(B)).
$$
Applying $(1\times \pi)_*$, this leads to a sheaf homomorphism,
$$
\partial G:\zeta_{\PP V}^* T_{\text{ev},\PP V} \rightarrow \text{pr}_B^*
\OO_B(d)\{t_0^j t_1^{d-1-j} | 0\leq j \leq d-1\}.
$$
By the deformation theory of stable maps, this $2$-term complex is
quasi-isomorphic to the good $[0,1]$ truncation of
$\mathbb{R}\textit{Hom}$ of the pullback by $\zeta_{\mb{X}}$ of the
relative cotangent complex of $\text{ev}:\Kbm{0,1}{\mb{X},1}
\rightarrow \mb{X}$.  In particular, the image of $\zeta_{\mb{X}}$ is
contained in $\Kbm{0,1}{\mb{X},1}_{\text{ev}}$ iff $\partial G$ is a
surjective sheaf homomorphism.  And in this case, $\zeta_{\mb{X}}^*
T_{\text{ev}, \mb{X}}$ is isomorphic to the kernel of $\partial G$,
cf. ~\cite[Rmk. 4.4]{HS2}.

\medskip\noindent
There is a canonical splitting $\OO_{D} \cong K \oplus
\mathfrak{m}\OO_D$.  The domain and target of $\partial G$ are
pullbacks under $\text{pr}_B$ of locally free sheaves on $B$.
Therefore $\partial G$ has a canonical splitting $\partial G =
\partial G_0 + \partial G_{\mathfrak{m}}$. 
By Nakayama's lemma, to
prove that $\partial G$ is surjective, it suffices to prove that
$\partial G_0$ is surjective after reducing modulo $\mathfrak{m}$.  
The homomorphism $\partial G$ is computed in Figure~\ref{table-1}.

\begin{figure}
$$
\begin{array}{|c|c|c|}
\hline
\text{Element} & \partial G_0 & \partial G_{\mathfrak{m}} \\
\hline
 & & \\
\mb{c} & 0 & \sum_{l} (2v^l \gamma_z + \sum_{(i,j)} u_{(i,j)}^l
\gamma_{(i,j)}) f^*X_l \\
\mb{b}_{(i,j)} & S_0^iS_1^{d-1-i}\cdot T_0^j T_1^{d-1-j} & - \\
\mb{a}_0 & S_1^{d-1}\cdot T_1^{d-1} & - \\
\mb{a}_1 & -S_1^{d-1}\cdot T_0T_1^{d-2} & - \\
\mb{a}_2 & -S_0S_1^{d-2} \cdot T_1^{d-1} & - \\
\mb{a}_3 & S_0S_1^{d-2}\cdot T_0T_1^{d-2} & - \\
\mb{e}_0 & S_1^{d-2}\cdot T_1^{d-2} & - \\
\mb{e}_1 & -S_1^{d-2}\cdot T_0T_1^{d-2} & - \\
\hline
\end{array}
$$ 
\caption{The homomorphism $\partial G$}
\label{table-1}
\end{figure}
\marpar{table-1}

\medskip\noindent
Some explanation of Figure~\ref{table-1} is in order.  The map being described
in the first 6 rows is the map on global sections induced by the sheaf
homomorphism $V\otimes_K \OO_{D\times B} \rightarrow
\text{pr}_B^* \OO_B(d-1)\{T_0^j T_1^{d-1-j} | 0\leq j \leq d-1\}$
obtained by composing $\partial G$ with the canonical surjection
$V\otimes_K \text{pr}_B^* \OO_B(1) \rightarrow \zeta_{\PP V}^*
T_{\text{ev},\PP V}$ described above, and then twisting by
$\OO_B(-1)$.  Each of these rows is obtained by taking the partial
derivative of $G$ with respect to the corresponding dual coordinate,
and the pulling back by $f^*$.  

\medskip\noindent
The last two rows are the images of $\mb{e}_0$ and
$\mb{e}_1$ under the map on global section obtained from $\partial G$
by twisting by $\OO_B(-2)$.  The images are uniquely determined by the
images of the $\mb{a}_l$s together with the relation between
$\mb{e}_0$, $\mb{e}_1$ and the $\mb{a}_l$s.  

\medskip\noindent  
The images of $\mb{b}_{(i,j)}$ and $\mb{e}_0$, $\mb{e}_1$ under
$\partial G_{\mathfrak{m}}$ are not given because they will not be used.
The image of $\mb{c}$ will be used.  

\begin{prop} \label{prop-dGsurj}
\marpar{prop-dGsurj}
\begin{enumerate}
\item[(i)]
The sheaf homomorphism $\partial G$ is surjective.  Moreover, the
sheaf homomorphism obtained by twisting by $\text{pr}_B^*\OO_B(-1)$ is
surjective 
on global sections.  
\item[(ii)]
The image of $\zeta_{\mb{X}}$ is
contained in $\Kbm{0,1}{\mb{X},1}_{\text{ev}}$, and $\zeta_{\mb{X}}^*
T_{\text{ev}, \mb{X}}$ is generated by global sections.
\item[(iii)]
The morphism $\zeta_{\mb{X}}:D\times B \rightarrow \Kbm{0,1}{\mb{X},1}$ is a
minimal twisting family parametrized by $D$.
\item[(iv)]
For every $\ssf{s}\in B$, the restriction of $\zeta$ to $D_\ssf{s}$
defines a morphism from $D_\ssf{s}$ to $M_\ssf{s}$.
\end{enumerate}
\end{prop}

\begin{proof}
\textbf{(i):}
Because $\text{pr}_B^* \OO_B(d)\{ T_0^j T_1^{d-1-j} | 0\leq j \leq
d-1\}$ is generated by global sections even after twisting by
$\text{pr}_B^* \OO_B(-1)$, it suffices to prove the sheaf
homomorphism obtained from $\partial G$ after twisting by
$\text{pr}_B^*\OO_B(-1)$ 
is surjective on global sections.
Every monomial $S_0^i S_1^{d-1-i} T_0^j T_1^{d-1-j}$ in
$H^0(B,\OO_B(d-1))\{T_0^j T_1^{d-1-j}| 0 \leq j \leq d-1 \}$ 
occurs as the image under
$\partial G_0$ of either $\mb{b}_{(i,j)}$, if $(i,j)\in I_d$,
or of $S_0\mb{e}_0$, $S_1\mb{e}_0$, $S_0\mb{e}_1$ or $S_1\mb{e}_1$.

\medskip\noindent
\textbf{(ii):}
As discussed above, because $\partial G$ is surjective the image of
$\zeta_{\mb{X}}$ is contained in $\Kbm{0,1}{\mb{X},1}_{\text{ev}}$ and
$\zeta_{\mb{X}}^* T_{\text{ev},\mb{X}}$ is the kernel of $\partial G$.
Because 
$\partial G$ is surjective on global sections after twisting by
$\text{pr}_B^* \OO_B(-1)$, $H^1(D \times B, \zeta_{\mb{X}}^*
T_{\text{ev},\mb{X}})$ is a subspace of $H^1(D\times B,\zeta_{\PP V}^*
T_{\text{ev}, \PP V})$.  But this is zero.  Therefore
$\zeta_{\mb{X}}^* T_{\text{ev},\mb{X}}$ is generated by global
sections.

\medskip\noindent
\textbf{(iii):}
By (ii), the image of $\zeta_{\mb{X}}$ is contained in
$\Kbm{0,1}{\mb{X},1}_{\text{ev}}$, and the pullback of
$T_{\text{ev},\mb{X}}$ is generated by global sections.  By
construction, $f:D\times \Sigma \rightarrow D\times \mb{X}$ is a
closed immersion whose image is a family of smooth quadric surfaces.  Also
$f \circ (1,\sigma):D\times B \rightarrow D\times \mb{X}$ is a
closed immersion whose image is a family of lines.  

\medskip\noindent
The only thing left to check is that the image of $f\circ (1,\sigma):D\times B
\rightarrow D\times \mb{X}$ is a family of \emph{free} lines.  
The computation so far
is symmetric with respect to the involution that permutes
$\PP^1_s$ and $\PP^1_t$.  By (ii), every fiber of
$\text{pr}_{\PP^1_s}:\PP^1_s \times \PP^1_t \rightarrow \PP^1_s$ is
mapped under $f$ to a free line.  So also every fiber of
$\text{pr}_{\PP^1_t}$ is mapped to a free line.  Since the image of
  $\sigma$ is the fiber over $[1:0]\in \PP^1_t$, $f\circ
  (1,\sigma):D \times\PP^1_s \rightarrow D\times\mb{X}$ maps to a family
  of free lines.

\medskip\noindent
\textbf{(iv):}  This follows from (iii) and Lemma~\ref{lem-Ds}.
\end{proof}

\begin{notat} \label{notat-kernel}
\marpar{notat-kernel}
A kernel of $\partial G$ is,
$$
\iota: \text{pr}_B^* \OO_B\{ \mb{f}_{(i,j)} | (i,j)\in I_d, i\leq d-2 \}
\oplus \text{pr}_B^* \OO_B\{\mb{g}_0,\mb{g}_1 \} \oplus \text{pr}_B^*
\OO_B(1) \{ \mb{h} \} \rightarrow \zeta_{\PP V}^* T_{\text{ev},\PP V},
$$
defined by,
$$
\lt\{
\begin{array}{ccccc}
\mb{f}_{(i,j)} & \mapsto & S_0 \mb{b}_{(i,j)} - S_1 \mb{b}_{(i+1,j)} &
+ & \mb{f}_{(i,j),\mathfrak{m}} \\ 
\mb{g}_0 & \mapsto & S_0^2 \mb{e}_0 - S_1 \mb{b}_{(2,0)} & +
& \mb{g}_{0,\mathfrak{m}} \\
\mb{g}_1 & \mapsto & S_0^2 \mb{e}_1 + S_1 \mb{b}_{(2,1)} & +
& \mb{g}_{1,\mathfrak{m}} \\
\mb{h} & \mapsto & \mb{c} & + & \mb{h}_{\mathfrak{m}}
\end{array} \rt. ,
$$
for some choice of 
$$
\mb{f}_{(i,j),\mathfrak{m}},
\mb{g}_{l,\mathfrak{m}} \in \mathfrak{m}\cdot H^0(D\times B,
\zeta_{\PP V}^* T_{\text{ev},\PP V}), \text{ and }
\mb{h}_{\mathfrak{m}} \in \mathfrak{m}\cdot
H^0(D\times B, \zeta_{\PP V}^* T_{\text{ev},\PP V}\otimes
\text{pr}_B^* \OO_B(-1)).
$$
\end{notat}

\medskip\noindent
There are many choices of
$\mb{f}_{(i,j),\mathfrak{m}}$ and $\mb{g}_{l,\mathfrak{m}}$ such that
$\iota$ is an isomorphism to the kernel of $\partial G$.
The element $\mb{h}_{\mathfrak{m}}$ is more canonical, and this is
crucial to the computation.  The space of global sections,
$$
H^0(D \times B,\zeta_{\mb{X}}^* T_{\text{ev},\mb{X}}\otimes
\text{pr}_B^*\OO_B(-1)),
$$
is a free $\OO_D$-module of rank $1$.  So $\mb{h}$ is unique up
to multiplication by an element in $1 + \mathfrak{m}$.  Therefore
$\mb{h}_{\mathfrak{m}}$ is a well-defined global section of,
$$
\mathfrak{m} \otimes_K \lt( \zeta_{\PP V,0}^*
T_{\text{ev}, \PP V} \otimes \OO_B(-1) / \OO_{B}\{\mb{c}\}\rt) =
\mathfrak{m} \otimes_K \textit{Hom}_{\OO_B}(\OO_{B}(1)\{\mb{c}\},
\zeta_{\PP V,0}^* T_{\text{ev}, \PP V}/\OO_{B}(1)\{\mb{c}\}).
$$
A lift of $\mb{h}_{\mathfrak{m}}$ to $\mathfrak{m} \otimes_K 
\zeta_{\PP V,0}^* T_{\text{ev}, \PP V} \otimes \OO_B(-1)$ is any
global section,
$$
\sum_{(i,j)\in I_d} \mb{h}_{\mathfrak{m},(i,j)} \mb{b}_{(i,j)} +
\mb{h}_{\mathfrak{m},0} \mb{e}_0 + \mb{h}_{\mathfrak{m},1} \mb{e}_1,
$$
whose image under $\partial{G}_0$ equals
$-\partial{G}_{\mathfrak{m}}(\mb{c})$. 




\section{The derivative map} \label{sec-deriv}
\marpar{sec-deriv}

\noindent
In the last section, a canonical global section $\mb{h}_\mathfrak{m}$
was identified in,
$$
\mathfrak{m}\OO_{D} 
\otimes_K \textit{Hom}_{\OO_B}(\OO_{B}(1)\{\mb{c}\},
\zeta_{\PP V,0}^* T_{\text{ev}, \PP V}/\OO_{B}(1)\{\mb{c}\}).
$$
In this section, $\mb{h}_\mathfrak{m}$ is restricted to
$\mathfrak{m}\OO_{D_\ssf{s}}$ and it is shown that this gives rise to
the restriction to $T_0 D_\ssf{s}$ of the
derivative map $dq_\ssf{s}|_{\zeta_0}$; this restriction is denoted
$d'q_\ssf{s}$.  

\medskip\noindent
Fix $\ssf{s} \in B$.  
Restrict from $D$ to the closed subscheme
$D_\ssf{s}$.
By Lemma~\ref{lem-Ds},
$\zeta_{\mb{X}}(D_\ssf{s}\times \{ \ssf{s} \})$ equals
$\zeta_0(\ssf{s})$ equals $[(L_\ssf{s},x_\ssf{s})]$.
Therefore, the map $\iota|_{D_\ssf{s}}$ gives a map of
$\OO_{D_\ssf{s}}$-modules,
$$
\text{pr}_B^*\OO_B(1)\{\mb{h}\} \otimes_{\OO_B} \kappa(\ssf{s}) \rightarrow
T_{L_\ssf{s},x_\ssf{s}} =
H^0(L_\ssf{s},N_{L_\ssf{s}/\mb{X}}(-x_\ssf{s})).
$$  
This in turn induces a $K$-linear map,
$$
d'q_\ssf{s}: T_0 D_\ssf{s} = \text{Hom}_K(\mathfrak{m}\OO_{D_\ssf{s}},K)
\rightarrow \text{Hom}_K(\OO_B(1)\{\mb{h}\}\otimes \kappa(\ssf{s}),
T_{L_\ssf{s}, x_\ssf{s}} / \OO_B(1)\{\mb{h}\}\otimes \kappa(\ssf{s})  ).
$$

\medskip\noindent
The map $d'q_\ssf{s}$ has another description.  The element
$\mb{h}_{\mathfrak{m}}$ from the last section gives a global section
$\mb{h}_{\mathfrak{m},\ssf{s}}$ of,
$$
\mathfrak{m}\OO_{D_\ssf{s}} 
\otimes_K \textit{Hom}_{\OO_B}(\OO_{B}(1)\{\mb{c}\},
\zeta_{\PP V,0}^* T_{\text{ev}, \PP V}/\OO_{B}(1)\{\mb{c}\}).
$$
This global section restricts to give an element in,
$$
\mathfrak{m}\OO_{D_\ssf{s}}
\otimes_K \textit{Hom}_{\OO_B}(\OO_{B}(1)\{\mb{c}\},
\zeta_{\PP V,0}^* T_{\text{ev}, \PP V}/\OO_{B}(1)\{\mb{c}\})
\otimes_{\OO_B} \kappa(\ssf{s}).
$$
The image of this element in,
$$
\mathfrak{m}\OO_{D_\ssf{s}} \otimes_K \OO_{B}(d-1)\{T_0^j T_1^{d-1-j}|
0\leq j \leq d-1 \} \otimes_{\OO_B} \kappa(\ssf{s}),
$$
is zero.  Therefore this is an element in,
$$
\mathfrak{m}\OO_{D_\ssf{s}}
\otimes_K \textit{Hom}_{\OO_B}(\OO_{B}(1)\{\mb{h}\},
\zeta_{\mb{X},0}^* T_{\text{ev}, \mb{X}}/\OO_{B}(1)\{\mb{h}\})
\otimes_{\OO_B} \kappa(\ssf{s}).
$$
By adjointness of $\otimes$ and $\text{Hom}$, this induces the
homomorphism $d'q_\ssf{s}$ above.

\medskip\noindent
The map $d'q_\ssf{s}$ can be made more explicit.  Let $\theta$ be an
element in $T_0 D_\ssf{s} = \text{Hom}_K(\mathfrak{m}
\OO_{D_\ssf{s}},K)$.  For each $(i,j) \in I_d$, the element $\langle
\mb{h}_{\mathfrak{m}, (i,j)}, \theta \rangle$ is just an element in $K$.
For $l=0,1$, the element $\langle \mb{h}_{\mathfrak{m}, l}, \theta \rangle$
is an element of $H^0(B,\OO_B(1))$.  The image in
$H^0(B,\OO_B(1)\otimes \OO_B(-\ssf{s}))$ is just an element in $K$.
So $d'q_\ssf{s}(\theta)$ is really just a vector of elements in $K$.  
All of these elements of $K$ are uniquely determined by,
$$
- \langle \partial G_{\mathfrak{m}}(\mb{c}), \theta \rangle \in
H^0(B,\OO_B(d-1)\{ T_0^j T_1^{d-1-1} | 0\leq j \leq d-1\} ).
$$

\medskip\noindent
The algorithm for computing $d'q_\ssf{s}$ is the following.
\begin{enumerate}
\item[Step 1]  Choose an ordered basis for $\mathfrak{m}
  \OO_{D_\ssf{s}}$ and a dual basis for $T_0 D_\ssf{s}$.
\item[Step 2]  For each basis element $w$ of $\mathfrak{m}
  \OO_{D_\ssf{s}}$ with dual basis element $\theta$, compute $\langle
  \partial G_{\mathfrak{m}}(\mb{c}), \theta \rangle$, i.e., compute
  the coefficient of $w$ in the expression $\partial
  G_{\mathfrak{m}}(\mb{c})$ in the vector space,
  $\mathfrak{m}\OO_{D_\ssf{s}} \otimes H^0(B,\OO_B(d-1)\{ T_0^j
  T_1^{d-1-j} | 0\leq j \leq d-1 \})$.  
\item[Step 3]  Find the unique element in $H^0(B,\zeta_{\PP V,0}^*
  T_{\text{ev},\PP V}\otimes \OO_B(-1))$ that is mapped by $\partial
  G_0$ 
  to $-\langle
  \partial G_{\mathfrak{m}}(\mb{c}), \theta \rangle$, 
  i.e., compute all the coefficients
  of $\mb{b}_{(i,j)}$, and of $S_0\mb{e}_0$, $S_1\mb{e}_0$,
  $S_0\mb{e}_1$ and $S_1\mb{e}_1$.  
\item[Step 4]  Compute the image of this element in $\zeta_{\PP V,0}^*
  T_{\text{ev}, \PP V} \otimes \OO_B(-1) \otimes_{\OO_B}
  \kappa(\ssf{s})$.  In essence, consider the element modulo
  $H^0(B,\OO_B(1)\{\mb{e}_0, \mb{e}_1\}\otimes \OO_B(-\ssf{s}))$.  
\item[Step 5]  Compute the element of $\zeta_{\mb{X},0}^*
  T_{\text{ev},\mb{X}}\otimes \OO_B(-1)/\OO_B\{\mb{h}\}
  \otimes_{\OO_B} \kappa(\ssf{s})$ that is mapped by $\iota\otimes
  \kappa(\ssf{s})$ to this element, i.e.,
  compute all the coefficients of $\mb{f}_{(i,j)}$, $\mb{g}_0$ and
  $\mb{g}_1$.  
\end{enumerate}

\medskip\noindent
\begin{prop} \label{prop-1}
\marpar{prop-1}
There exist polynomials $\gamma_z$ and $\gamma_{(i,j)}$, for $(i,j)\in
I_d$, and  
there exists $\ssf{s}\in \PP^1_s$ such that
$q_\ssf{s}$ is smooth at $\zeta_{\mb{X}}|_0$.  
\end{prop}

\medskip\noindent
By Lemma~\ref{lem-q}, the domain and target of $q_\ssf{s}$ are
smooth.  So, by the Jacobian criterion, $q_\ssf{s}$ is smooth at
$\zeta_0$ iff
$dq_\ssf{s}|_{\zeta_0}$ is surjective.  To prove this is surjective,
it suffices to prove the restriction $d'q_\ssf{s}$ is surjective.

\section{Proof of the proposition} \label{sec-proof}
\marpar{sec-proof}

\noindent
The point $\ssf{s}$ is $[1:0]$.  The polynomials $\gamma_{(i,j)}$ and
$\gamma_z$ are given in Figure~\ref{table-2}.
The terms $C_{(i,j)}$, $C_j$, $C_{1,a}$, $C_{1,b}$, $C_{z,a}$ and
$C_{z,b}$ 
are elements of $K$ that are generic (this is made more precise
later).  
For any index
$(i,j)$ not appearing above, $C_{(i,j)}$ is defined to be $0$.
The remainder of this section carries out the algorithm of the
previous section.

\begin{figure}
$$
\begin{array}{|c|c|c|}
\hline
\text{Index} & \text{Polynomial} & f^*\text{Polynomial} \\
\hline
\gamma_{(i,j)}, &  & \\
(i,j) \in I_d, i,j\leq d-2 & C_{(i,j)} X_0^k
X_1^{i-k}X_2^{j-k}X_3^{d-2-i-j+k} & C_{(i,j)}S_0^iS_1^{d-2-i}\cdot
T_0^j T_1^{d-2-j} \\ 
\hline
\gamma_{(d-1,j)}, &  & \\
2\leq j \leq d-1 &
C_j X_1^{j-1} X_3^{d-1-j} & C_j S_0^{j-1} S_1^{d-1-j} T_1^{d-2} \\
\hline
 & & \\
\gamma_{(d-1,1)} & C_{1,a} X_1X_3^{d-3} + C_{1,b} X_2 X_3^{d-3} &
C_{1,a} S_0S_1^{d-3}\cdot T_1^{d-2} + C_{1,b}S_1^{d-2}\cdot T_0T_1^{d-3}  \\
\hline
 & & \\
\gamma_z & C_{z,a} X_0X_3^{d-3} + C_{z,b} X_1X_3^{d-3} &
C_{z,a} S_0S_1^{d-3}\cdot T_0T_1^{d-3} + C_{z,b} S_0S_1^{d-3}\cdot T_1^{d-2} \\
\hline
\gamma_{(i,d-1)}, & & \\
0\leq i \leq d-2 & 0 & 0 \\
\hline
\gamma_{(d-2,d-1)} & 0 & 0\\
\hline
\gamma_{(d-1,0)} & 0 & 0 \\
\hline
\end{array}
$$
\caption{The polynomials $\gamma$}
\label{table-2}
\end{figure}
\marpar{table-2}

\medskip\noindent
\textbf{Step 1:}
Define $J_d = \{(i,j)\in I_d| (i,j) \neq (i,0), (0,2),(1,2)\}$.  
There is an isomorphism of $K$-vector spaces,
$$
\mathfrak{m}\OO_{D_\ssf{s}} \rightarrow K\{v^2,v^3\}\cup \{w_{(i,j)} | (i,j)
\in J_d\},
$$
given in Figure~\ref{table-3}.

\begin{figure}
$$
\begin{array}{|c|c|c|}
\hline
\text{Variable} & \text{Image} & \text{Range} \\
\hline & & \\
u_{(i,j)}^0 & 0 &  - \\
u_{(i,j)}^1 & 0 &  - \\
u_{(i,j)}^2 & -w_{(i,j+1)} & j\leq d-2 \\
u_{(i,d-1)}^2 & 0 & -  \\
u_{(i,j)}^3 & w_{(i,j)} & (i,j) \in J_d \\
u_{(i,j)}^3 & 0 & (i,j) = (i,0), (0,2), (1,2) \\
v^0 & 0 & - \\
v^1 & 0 & - \\
v^2 & v^2 & - \\
v^3 & v^3 & - \\
\hline
\end{array} 
$$
\caption{Basis for $\mathfrak{m}\OO_{D_s}$}
\label{table-3}
\end{figure}
\marpar{table-3}

\medskip\noindent
\textbf{Step 2:}
With respect to the isomorphism above, $\mb{h}_{\mathfrak{m}}$ equals,
$$
\mb{h}_{\mathfrak{m}} = \sum_{(i,j)\in J_d} w_{(i,j)}\cdot
\langle \mb{h}_{\mathfrak{m}}, w_{(i,j)}^\vee \rangle + 
v^1\cdot \langle \mb{h}_{\mathfrak{m}}, (v^1)^\vee \rangle +
v^2\cdot \langle \mb{h}_{\mathfrak{m}}, (v^2)^\vee \rangle
$$
given in Figure~\ref{table-4}.  Each computation in
Figure~\ref{table-4} is straightforward, and is left to the reader.

\begin{figure}
$$
\begin{array}{|c|c|}
\hline
\text{Element } w & \langle \mb{h}_{\mathfrak{m}}, w^\vee \rangle \\
\hline 
w_{(i,j)}, & \\
(i,j)\in J_d, i\leq d-2 & (C_{(i,j)} -
C_{(i,j-1)})S_0^iS_1^{d-1-i}T_0^jT_1^{d-1-j} \\
\hline
w_{(d-1,j)}, &  \\
3\leq j \leq d-1 &  C_j S_0^{j-1}S_1^{d-j}T_1^{d-1} -
C_{j-1}S_0^{j-2}S_1^{d+1-j}T_0T_1^{d-2} \\
\hline &  \\
w_{(d-1,2)} & -C_{1,a} S_0S_1^{d-2}T_0T_1^{d-2} -C_{1,b}
S_1^{d-1}T_0^2T_1^{d-3} + C_2 S_0S_1^{d-2}T_1^{d-1} \\
\hline &  \\
w_{(d-1,1)} & C_{1,a}S_0S_1^{d-2}T_1^{d-1} +
C_{1,b}S_1^{d-1}T_0T_1^{d-2} \\
\hline &  \\
v^2 & C_{z,a} S_0S_1^{d-2}T_0^2T_1^{d-3}
+C_{z,b}S_0S_1^{d-2}T_0T_1^{d-2} \\
\hline &  \\
v^3 & C_{z,a} S_0S_1^{d-2}T_0T_1^{d-2} + C_{z,b} S_0S_1^{d-2}T_1^{d-1} \\
\hline 
\end{array}
$$ 
\caption{Computation of $\langle \mb{h}_{\mathfrak{m}}, w^\vee \rangle$}
\label{table-4}
\end{figure}
\marpar{table-4}

\medskip\noindent
\textbf{Step 3:}  For each $w$ in the basis of
$\mathfrak{m}\OO_{D_\ssf{s}}$, an element in $H^0(B,\zeta_{\PP V,0}^*
T_{\text{ev},\PP V}\otimes \OO_B(-1))$ that maps under $\partial G_0$
to $-\langle \mb{h}_{\mathfrak{m}}, w^\vee \rangle$ is given in 
Figure~\ref{table-5}

\begin{figure}
$$
\begin{array}{|c|c|}
\hline
\text{Element } w & - \partial{G}_0^{-1}
\langle \mb{h}_{\mathfrak{m}}, w^\vee \rangle \\
\hline 
w_{(i,j)}, &  \\
(i,j)\in J_d, i\leq d-2 & (-C_{(i,j)}+C_{(i,j-1)}) \mb{b}_{(i,j)}  \\
\hline
w_{(d-1,j)}, & \\
4\leq j \leq d-1 &  -C_j \mb{b}_{(j-1,0)} + C_{j-1} \mb{b}_{(j-2,1)}  \\
\hline & \\
w_{(d-1,3)} & - C_3 \mb{b}_{(2,0)}  - C_2 S_0 \mb{e}_1 \\
\hline & \\
w_{(d-1,2)} & -C_{1,a}S_0 \mb{e}_1 + C_{1,b} \mb{b}_{(0,2)} - C_2 S_0
\mb{e}_0 \\
\hline & \\
w_{(d-1,1)} & -C_{1,a} S_0 \mb{e}_0 + C_{1,b} S_1 \mb{e}_1 \\
\hline & \\
v^2 & -C_{z,a} \mb{b}_{(1,2)} + C_{z,b} S_0 \mb{e}_1 \\
\hline & \\
v^3 & C_{z,a}S_0 \mb{e}_1 - C_{z,b} S_0 \mb{e}_0 \\
\hline 
\end{array}
$$
\caption{Computation of $-\partial{G}_0^{-1} \langle
  \mb{h}_{\mathfrak{m}}, w^\vee \rangle$}
\label{table-5}
\end{figure}
\marpar{table-5}

\medskip\noindent
\textbf{Step 4:}  The images of these elements in $\zeta{\PP V,0}^*
T_{\text{ev}, \PP V}\otimes \OO_B(-1) \otimes_{\OO_B} \kappa(\ssf{s})$
are obtained by setting $S_1$ equal to $0$.  This only effects the
term for $w_{(d-1,1)}$.

\medskip\noindent
\textbf{Step 5:}
The vector space map,
$$
\zeta_{\mb{X},0}^* T_{\text{ev},\mb{X}} \otimes
  \OO_B(-1)/\OO_B\{\mb{h}\} \otimes_{\OO_B} \kappa(\ssf{s})
  \rightarrow
\zeta_{\PP V,0}^* T_{\text{ev},\PP V} \otimes
  \OO_B(-1)/\OO_B\{\mb{c}\} \otimes_{\OO_B} \kappa(\ssf{s}),
$$
is given by,
$$
\lt\{
\begin{array}{ccc}
\mb{f}_{(i,j)} & \mapsto & S_0 \mb{b}_{(i,j)} \\
\mb{g}_0 & \mapsto & S_0^2 \mb{e}_0 \\
\mb{g}_1 & \mapsto & S_0^2 \mb{e}_1 \\
\end{array} \rt.
$$
For each basis element $w \in \mathfrak{m}
\OO_{D_\ssf{s}}$, $d'q_\ssf{s}(w^\vee)$ is the preimage of the element
in the previous table.  This is given in Figure~\ref{table-6}.

\begin{figure}
$$
\begin{array}{|c|c|}
\hline
\text{Element } w & d'q_\ssf{s}(w^\vee) \\
\hline 
w_{(i,j)}, & \\
(i,j)\in J_d, i\leq d-2 &
(-C_{(i,j)} + C_{(i,j-1)})
\frac{1}{S_0} \mb{f}_{(i,j)} \\
\hline
w_{(d-1,j)}, & \\
4\leq j \leq d-1 &  
-C_j \frac{1}{S_0} \mb{f}_{(j-1,0)} +
C_{j-1}\frac{1}{S_0} \mb{f}_{(j-2,1)} \\
\hline & \\
w_{(d-1,3)} & 
-C_3 \frac{1}{S_0} \mb{f}_{(2,0)} - C_2
\frac{1}{S_0}  \mb{g}_1 \\
\hline & \\
w_{(d-1,2)} & 
-C_{1,a} \frac{1}{S_0} \mb{g}_1 + C_{1,b}
\frac{1}{S_0} \mb{f}_{(0,2)} 
- C_2 \frac{1}{S_0}
 \mb{g}_0 \\  
\hline & \\
w_{(d-1,1)} &
-C_{1,a} \frac{1}{S_0} \mb{g}_0 \\
\hline & \\
v^2 & 
-C_{z,a}
\frac{1}{S_0} \mb{f}_{(1,2)} + C_{z,b}
\frac{1}{S_0} \mb{g}_1 \\
\hline & \\
v^3 & 
C_{z,a}
\frac{1}{S_0}  \mb{g}_1 - C_{z,b} \frac{1}{S_0}
 \mb{g}_0 \\  
\hline 
\end{array}
$$
\caption{Computation of $d'q_{\ssf{s}}$}
\label{table-6}
\end{figure}
\marpar{table-6}

\medskip\noindent
From the first row, the image of $d'q_{\ssf{s}}$ contains
$\frac{1}{S_0} \mb{f}_{(i,j)}$ for all $(i,j)\in
J_d$ and $i \leq d-2$.  In particular, $(j-2,1)$ is such an index for
$4\leq j \leq d-1$.  So, from the second row, the term
$-C_{j-1} \frac{1}{S_0} \mb{f}_{(j-2,1)}$ is
already in the image of $d'q_{\ssf{s}}$.  Therefore $\frac{1}{S_0}
\mb{f}_{(j-1,0)}$ is in the image of
$d'q_{\ssf{s}}$ for $4\leq j \leq d-1$.  

\medskip\noindent
Among the terms $\mb{f}_{(i,j)}$,
the only ones not covered by the first 2 rows are $(i,j)=(0,2),(1,2)$
and $(2,0)$.

\medskip\noindent
From the fifth row, $\frac{1}{S_0}\mb{g}_0$ is in
the image of $d'q_{\ssf{s}}$.  Together with the last row, then
$\frac{1}{S_0} \mb{g}_1$ is in the image of
$d'q_{\ssf{s}}$.  With the fourth row, 
$\frac{1}{S_0} \mb{f}_{(0,2)}$ is in the image.  With the sixth row,
$\frac{1}{S_0} \mb{f}_{(1,2)}$ is in the image.
And with the third row, $\frac{1}{S_0} \mb{f}_{(2,0)}$ is in the
image.  Therefore $d'q_{\ssf{s}}$ is surjective. 

\section{A spanning minimal twisting family} \label{sec-pfmain}
\marpar{sec-pfmain}

\noindent
If $d=1$ or $2$, any minimal twisting family is already very twisting.

\begin{prop} \label{prop-twistlines}
\marpar{prop-twistlines}
Let $\mb{X}\subset \PP^n$ be a general hypersurface of degree $d\geq
3$.  If $n\geq d^2$, there exists a point $\ssf{s}\in \PP^1$, a point
$[L,x]\in \Kbm{0,1}{\mb{X},1}$ and a minimal twisting morphism
$\zeta_{\mb{X},0}:(\PP^1,\ssf{s}) \rightarrow
(\Kbm{0,1}{\mb{X},1},[L,x])$ such that $q_{\ssf{s},[L,x]}$ is spanning
at $\zeta_{\mb{X},0}$.  
\end{prop}

\medskip\noindent
Let $d\geq 3$ and let $n'\geq d^2$.  Denote $n=d^2$.  Choose a linear
subspace $\PP^{n} \subset \PP^{n'}$.  Let $\mb{X} \subset \PP^{n}$ be
the hypersurface from Section~\ref{sec-proof}.  Let $\mb{X}' \subset
\PP^n$ be a cone over $\mb{X}$.  Let $\ssf{s}=[1:0]$ and let
$\zeta_{\mb{X}}:D_\ssf{s} \times B \rightarrow \Kbm{0,1}{\mb{X},1}$
be the minimal twisting family from Section~\ref{sec-not}.  This is
also a minimal twisting family for $\mb{X}'$; denote by
$\zeta_{\mb{X}'}: D_\ssf{s} \times B \rightarrow \Kbm{0,1}{\mb{X}',1}$
the corresponding morphism.  There is a short exact sequence,
$$
\begin{CD}
0 @>>> \zeta_{\mb{X}}^* T_{\text{ev},\mb{X}} @>>> \zeta_{\mb{X}'}^*
T_{\text{ev}, \mb{X}'} @>>> \text{pr}_B^* \OO_B(1)^{n'-n} @>>> 0.
\end{CD}
$$
Therefore, the quotient of $T_{\text{ev},\mb{X'}} \otimes
\kappa([L_\ssf{s},x_\ssf{s}])$ by $\text{pr}_B^* \OO_B(1)^{n'+1-d^2}$
is canonically isomorphic to the quotient of $T_{\text{ev},\mb{X}}
\otimes \kappa([L_\ssf{s},x_\ssf{s}])$ by $\text{pr}_B^*
\OO_B(1)^{n+1-d^2}$.  By Proposition~\ref{prop-1}, 
$d'q_\ssf{s}$ surjects to this quotient.  Therefore the morphism
$q_\ssf{s}$ for $\mb{X}'$ is smooth at $\zeta_0$; i.e., $q_\ssf{s}$ is
spanning.  By Proposition~\ref{prop-q}, for a general hypersurface of
degree $d$ in $\PP^{n'}$, there exists a very twisting family of
lines.  

\section{Can we do better?} \label{sec-misc}
\marpar{sec-misc}

\noindent
What is the weakest inequality such that the spaces
$\Kbm{0,0}{\mb{X},e}$ are rationally connected?  What is the weakest
inequality such that there exists a very twisting 1-morphism
$\zeta:\PP^1 \rightarrow \Kbm{0,1}{\mb{X},1}$?  The weakest inequality
such that there exists a very twisting morphism is $n\geq d^2$,
cf. Proposition~\ref{prop-ineq}.  But rational connectedness holds
under slightly weaker hypotheses.  The basic idea is contained in the
following result.

\begin{thm}[de Jong] \label{thm-uniruled}
\marpar{thm-uniruled}
Let $K$ be an algebraically closed field with $\text{char}(K)=0$.
Let $(d,n)$ be positive integers such that either $d=1$ and $n\geq 2$  
or $d\geq 2$ and $n+1 \geq d^2$.
Let $\mb{X} \subset \PP^n$ be a general hypersurface of degree $d$.
For every $e\geq 1$, every irreducible component of
$\Kbm{0,0}{\mb{X},e}$ is uniruled (or is a point if $(d,n,e)=(1,2,1)$).  
\end{thm}

\begin{proof}
If $d=1$ or $d=2$, this is easy and follows from a stronger
result of Kim and Pandharipande, cf. \cite{KP}.  In those cases 
there is an action of $\text{SL}_n$, resp. $\text{SO}_{n+1}$, 
on $\Kbm{0,0}{\mb{X},e}$.  With the one exception of $(d,n,e)=(1,2,1)$,
the stabilizer of a general point is a subgroup of positive codimension.
Because $\text{SL}_n$ and
$\text{SO}_{n+1}$ are rational, the orbit of a general point is
a positive dimensional, unirational variety; in particular it is
uniruled.  Because the orbit of a
general point is uniruled,
$\Kbm{0,0}{\mb{X},e}$ is uniruled.

\medskip\noindent
Let $d\geq 3$.  By ~\cite[Prop. 6.5]{HS2}, 
there exists a twisting family $\zeta:\PP^1 \rightarrow
\Kbm{0,1}{\mb{X},1}$ such that 
$\text{ev}\circ \zeta:\PP^1 \rightarrow \mb{X}$ is a line.  Let
$g:\PP^1 \rightarrow \PP^1$ be a finite morphism of degree $e$.  Then
$\zeta\circ g: \PP^1 \rightarrow \Kbm{0,1}{\mb{X},1}$ is also
twisting; i.e., there exists a degree $e$ stable map $\text{ev}\circ
\zeta\circ g: \PP^1 \rightarrow X$ that is \emph{twistable},
cf. ~\cite[Def. 4.7]{HS2}.  By ~\cite[Prop. 4.8]{HS2}, there is a
nonempty open subset of $\Kbm{0,0}{\mb{X},e}$ of twistable stable
maps.  By ~\cite{HRS2}, $\Kbm{0,0}{\mb{X},e}$ is irreducible,
therefore this open set is dense.  If $h:\PP^1 \rightarrow \mb{X}$ is
in this open set, then there exists a twisting morphism $\xi:\PP^1
\rightarrow \Kbm{0,1}{\mb{X},1}$ such that $\text{ev}\circ \xi:\PP^1
\rightarrow \mb{X}$ equals $h$.  Let $(\pi:\Sigma \rightarrow \PP^1,
\sigma:\PP^1 \rightarrow \Sigma, g:\Sigma \rightarrow \mb{X})$ be the
family inducing $\zeta$.  Then $\sigma(B) \subset \Sigma$ deforms in
its linear equivalence class.  This deformation defines a
non-constant, rational transformation $\PP^1 \dashrightarrow
\Kbm{0,0}{\mb{X},1}$ whose image contains $[h]$.  Therefore
$\Kbm{0,0}{\mb{X},1}$ is uniruled.
\end{proof}

\medskip\noindent
A. J. de Jong's proof is different, and proves a stronger result: for
\emph{every} smooth hypersurface, \emph{every} irreducible component of
$\Kbm{0,0}{\mb{X},e}$ is uniruled.  
A. J. de Jong's theorem motivated a re-investigation of
~\cite[Lem. 7.4]{HS2} and the proof of Theorem~\ref{thm-main}.

\medskip\noindent
In a forthcoming paper, by elaborating on the proof of 
Theorem ~\ref{thm-uniruled}, de Jong and I prove that if $n\geq
d^2+1$, then the spaces of curves are rationally connected.  However,
that proof does not give the existence of very twisting families of
lines.  

\medskip\noindent
This may seem inconsequential, but it is actually 
important for another purpose: de Jong and I have a method of
generalizing the Tsen-Lang 
theorem.  In that method, existence of a very twisting family of lines
plays an important role.  As the Tsen-Lang theorem holds for $n\geq
d^2$ but \emph{does not} hold for $n=d^2+1$, non-existence of a very
twisting family of lines appears to be signficant.  

\medskip\noindent
In the other direction, 
there is reason to expect that $\Kbm{0,0}{\mb{X},e}$ is not uniruled
if $n+2 \leq d^2$.

\begin{prop}~\cite{S1}  \label{prop-big}
\marpar{prop-big}
Let $K$ be an algebraically closed field with $\text{char}(K) = 0$.
Let $(d,n)$ be positive integers such that
$d < \min(n-3,\frac{n+1}{2})$ and $n+2 \leq d^2$ (if $n\geq 6$, the
conditions are $d\leq \frac{n}{2}$ and $n+2 \leq d^2$). 
Let $\mb{X}
\subset \PP^n$ be a general hypersurface of degree $d$.  For every $e \gg
0$ the canonical divisor class on $\Kbm{0,0}{\mb{X},e}$ is
\emph{big}.  If $e\leq n-d$, then the coarse moduli space
$\kbm{0,0}{\mb{X},e}$ is of general 
type (the singularities are \emph{canonical}).
\end{prop}

\bibliography{my}
\bibliographystyle{alpha}

\end{document}